\documentclass{infor}

\volume{0}
\issue{0}
\pubyear{2024}
\articletype{research-article}
\doi{0000}

\RequirePackage{fix-cm}
\usepackage{soul}
\usepackage{amsfonts}
\usepackage{graphicx}
\usepackage{graphics}
\usepackage{color}
\usepackage{graphicx}
\usepackage{rotating}
\usepackage{amsmath,amssymb}
\usepackage[ruled,vlined]{algorithm2e}
\usepackage{pdflscape}
\usepackage{multirow}
\usepackage{graphics}

\newcommand\mybox{\hbox to 0pt{}\hfill$\rlap{$\sqcap$}\sqcup$}
\newtheorem{lemma}{Lemma}

\theoremstyle{remark}

\begin{document}
\begin{frontmatter}
\pretitle{Research Article}

\title{A general framework for providing interval representations of Pareto optimal outcomes for large-scale bi- and tri-criteria MIP problems}

\runtitle{Providing interval representations of Pareto optimal outcomes}

\author[a]{\inits{G.}\fnms{Grzegorz} \snm{Filcek}
	\ead[label=e1]{grzegorz.filcek@pwr.edu.pl}\bio{bio1}}

\author[b]{\inits{J.}\fnms{Janusz} \snm{Miroforidis}\thanksref{f1}\ead[label=e2]{janusz.miroforidis@ibspan.waw.pl}\bio{bio2}}
\thankstext[type=corresp,id=f1]{Corresponding author.}
\address[a]{Ul. Wybrze\.ze Wyspia\'nskiego 27, 50-370 Wroc\l{}aw, \institution{Wroc\l{}aw University of Science and Technology}, \cny{Poland}}
\address[b]{Ul. Newelska 6, 01-447 Warszawa, \institution{Systems Research Institute, Polish Academy of Sciences}, \cny{Poland}}


\begin{abstract}
The Multi-Objective Mixed-Integer Programming (MOMIP) problem is one of the most challenging. To derive its Pareto optimal solutions one can use the well-known Chebyshev scalarization and  Mixed-Integer Programming (MIP) solvers. However, for a large-scale instance of the MOMIP problem, its scalarization may not be solved to optimality, even by state-of-the-art optimization packages, within the time limit imposed on the optimization.
If a MIP solver cannot derive the optimal solution within the assumed time limit, it provides the optimality gap, which gauges the quality of the approximate solution. However, for the MOMIP case, no information is provided on the lower and upper bounds of the components of the Pareto optimal outcome.
For the MOMIP problem with two and three objective functions, an algorithm is proposed to provide the so-called interval representation of the Pareto optimal outcome designated by the weighting vector when there is a time limit on solving the Chebyshev scalarization. Such interval representations can be used to navigate on the Pareto front. 
The results of several numerical experiments on selected large-scale instances of the multi-objective, multidimensional 0-1 knapsack problem illustrate the proposed approach. The limitations and possible enhancements of the proposed method are also discussed.

\end{abstract}


\end{frontmatter}

\section{Introduction}
\label{section1}
The derivation of optimal solutions to large-scale instances of the Mixed-Integer Programming (MIP) problem can be impossible within a reasonable time limit even for contemporary commercial MIP solvers, e.g., GUROBI (\cite{GUROBI}), CPLEX (\cite{CPLEX}). In this case, a MIP solver provides the optimality gap (MIP gap) that gauges the quality of the approximate solution, i.e., the last feasible solution (incumbent). This optimality gap is calculated based on the incumbent and the so-called MIP best bound.

In the case of the Multi-Objective MIP (MOMIP) problem, scalarization techniques and MIP solvers can be used to derive Pareto optimal solutions (see, e.g., \cite{Miettinen_1999}, \cite{Ehrgott_2005}).
Examples of applying MIP packages to solve multi-criteria decision problems are shown in, e.g., \cite{Ahmadi_2012}, \cite{Delorme_2014}, \cite{Eiselt_2014}, \cite{Oke_2015}, \cite{Samanlioglu_2013}.
As a scalarization technique, one can use the Chebyshev scalarization that guarantees the derivation of each (properly) Pareto optimal solution (see, e.g., \cite{Kaliszewski_2006}). Other advantages of using this scalarization in the context of decision-making and expressing the decision maker's preferences are discussed in, e.g., \cite{Miroforidis_2021_Markowitz}.

In the current work, we say that an instance of the MOMIP problem is large-scale if its Chebyshev scalarization cannot be solved to optimality by a MIP solver within an assumed time limit that is reasonable in the decision-making process.
The existence of this limit is justified in solving practical multi-criteria decision-making problems.
When there is a time limit on deriving a single Pareto optimal solution, the Chebyshev scalarization of the instance may not be solved to optimality.
The decision maker (DM) then obtains the incumbent, i.e. the approximation of the Pareto optimal solution, as well as the MIP gap of the single-objective optimization problem.
However, based on this information, the quality of the approximation of a single component (namely its lower and upper bounds) of the Pareto optimal outcome, i.e. the image of the Pareto optimal solution in the objective space, cannot be shown to the DM. And it is based on these components that the DM navigates on the Pareto front (set of Pareto optimal outcomes).
Fortunately, there is a method to provide the DM with such lower and upper bounds in the literature. 

In \cite{Kaliszewski_Miroforidis_2019_genbounds}, a general methodology for multi-objective optimization to provide lower and upper bounds on objective function values of a Pareto optimal solution designated by a vector of weights of the Chebyshev scalarization of a multi-objective optimization problem has been proposed.
The bounds form the so-called \emph{interval representation} of the Pareto optimal outcome. The DM can use interval representations instead of (unknown) Pareto optimal solutions,
to navigate on the Pareto front. To derive them, one needs the so-called \emph{lower shells} and \emph{upper shells} whose images in the objective space are finite two-sided approximations of the Pareto front
(see, e.g., \cite{Kaliszewski_Miroforidis_2014_two_sided}).

In \cite{Kal_Mir_2022_Probing}, it has been shown how to provide lower and upper shells to large-scale instances of the MOMIP problem. In that work, lower shells are composed of incumbents to the Chebyshev scalarization of the MOMIP problem derived within the time limit, and upper shells consist of elements that are solutions to the Chebyshev scalarization of a relaxation of the MOMIP problem.

However, there is a lack of an algorithmic method for deriving an upper shell that is necessary to calculate the interval representation of the Pareto optimal outcome designated by a given vector of weights of the Chebyshev scalarizing function. 
The idea of how to derive such useful upper shells for the MOMIP problem with two objective functions has been shown in our earlier works \cite{Kal_Mir_2021_Cooperative} and \cite{Miroforidis_2021_Markowitz}.

In the current work, we combine ideas from works \cite{Kaliszewski_Miroforidis_2019_genbounds}, \cite{Kal_Mir_2022_Probing}, \cite{Kal_Mir_2021_Cooperative}, and \cite{Miroforidis_2021_Markowitz}. 
For this reason, our work is an incremental one.
For the MOMIP problem with up to three objective functions, we propose an algorithmic method of deriving upper shells that can be used to calculate the interval representation of a single Pareto optimal outcome designated by a given vector of weights of the Chebyshev scalarizing function.
This opens the way for providing the DM with this representation when there is a time limit for deriving a single Pareto optimal solution.
Because of the need to derive the appropriate upper shells, additional time is needed for optimization, but as we show in numerical experiments, this time can be a fraction of the assumed time limit. To illustrate our method, we present results of several numerical experiments with selected large-scale instances of the multi-objective multidimensional 0-1 knapsack problem.

To our best knowledge, the method we propose is the only algorithmic method for determining the interval representation of the Pareto optimal outcome given by weights of the Chebyshev scalarizing function for large-scale instances of the MOMIP problem, assuming the existence of a time budget for optimization.
This method which is a generic framework for providing these interval representations is the main contribution of this work.

The current work is organized as following. In Section \ref{section2}, we formulate the MOMIP problem and we recall a method for the derivation of Pareto optimal solutions with the use of the Chebyshev scalarization. In Section \ref{bounds}, we briefly recall the theory of parametric lower and upper bounds. There, we also introduce the concept of the interval representation of the implicit Pareto optimal outcome as well as an indicator measuring its quality. 
In Section \ref{deriv_interval_repr}, we present two versions of an algorithm for deriving interval representations
of implicit Pareto optimal outcomes.
In Section \ref{numerical_experiments}, we conduct extensive numerical experiments as well as discuss their results.
In Section \ref{limitations}, we show the limitations of the proposed method as well as discuss how to eliminate them. Section \ref{final_remarks} contains some final remarks.

\section{Background}
\label{section2}

Let $X := {\mathbb{R}^{n_1}}\times{\mathbb {Z}}^{n_2}$, $n_1+n_2=n, n_2>0$, $x \in X$ denote a solution, $X_0:=\{x \in X\, |\, g_{p}(x) \leq b_{p},\, b_{p}\in \mathbb{R} \}$  the set of feasible solutions, where  $g_p : {\mathbb{R}}^n \rightarrow {\mathbb{R}}, \ p=1,\ldots,m, \ m \geq 1$. The MOMIP problem is defined as follows:
\begin{equation}
	\label{eq2.1}	
	\begin{array}{ll}
		\text{vmax} & f(x) \vspace{1mm}\\
		\text{s.t.} & x \in X_0 \,,\\
	\end{array}				
\end{equation}
where $f : {\mathbb{R}}^n \rightarrow {\mathbb{R}}^k$, $f = (f_1,\ldots,f_k)$, $f_l : {\mathbb{R}}^n \rightarrow {\mathbb{R}}, \ l=1,\ldots,k, \ k \geq 2$, are objective functions, and  "vmax" is the operator of deriving set $N$ that contains all Pareto optimal solutions in $X_0$. 
The set ${\mathbb{R}}^k$ is called the objective space. 
Solution $\bar{x}\in X_0$ is Pareto optimal, if for any $x\in X_0$, $f_l(x) \geq f_l(\bar{x}), \ l = 1,\ldots,k$, implies $f(x) = f(\bar{x})$. 
If $f_l(x) \geq f_l(\bar{x}), \ l = 1,\ldots,k$, and $f(x) \not= f(\bar{x})$, then $x$ dominates $\bar{x}$ ($\bar{x}$ is dominated) which is denoted by the relation $x \succ \bar{x}$. We say that element $f(x)$,  $x \in X_{0}$, is the outcome of $x$. Set $f(N)$ is called the Pareto front.

According to well-established knowledge (\cite{Ehrgott_2005}, \cite{Kaliszewski_2006}, \cite{KMP_2016_Springer}, \cite{Miettinen_1999}), the solution $x$ is Pareto optimal (actually, $x$ is properly Pareto optimal, see, e.g., \cite{Ehrgott_2005}, \cite{Kaliszewski_2006}, \cite{KMP_2016_Springer}, \cite{Miettinen_1999}) if and only if it solves the Chebyshev scalarization of problem (\ref{eq2.1}), namely

\begin{equation}
	\label{eq2.2}	
	\begin{array}{ll}
		\min & \max_l \lambda_l (y^*_l - f_l(x)) + \rho e^k(y^* - f(x)) \vspace{1mm}\\ \text{s.t.} & x \in X_0 \,,\\
	\end{array}				
\end{equation}
where weights $\lambda_l > 0, \ l= 1,\ldots,k$, \ $e^k = (1,1,\ldots,1)$, \ $y^*_l = \hat{y}_l + \varepsilon$, $\hat{y}_l = \max_{x \in X_0} f_l(x) < \infty, \ l=1,\ldots,k$, 
$\varepsilon > 0$, and $\rho$  is a positive "sufficiently small" number.
The linearized version of problem (\ref{eq2.2}) is the following.
\begin{equation}
\label{eq2.3}	
\begin{array}{ll}
	\min & s \vspace{1mm}\\
	\text{s.t.} & s \geq  \lambda_l (y^*_l - f_l(x)) + \rho e^k(y^* - f(x)), \ \ l=1,\ldots,k \,, \vspace{1mm}\\
	& x \in X_0 \,.
\end{array}				
\end{equation}
In the following, we will assume that Pareto optimal solutions come from solving problem (\ref{eq2.3}) with varying $\lambda=(\lambda_1, ..., \lambda_k)$. 

Given $\lambda$, $x^{P_{opt}}(\lambda)$ denotes the \emph{implicit Pareto optimal solution} designated
by $\lambda$ that is a solution which would be derived if problem (\ref{eq2.2}) with $\lambda$ were solved to optimality. $f(x^{P_{opt}}(\lambda))$ denotes the \emph{implicit Pareto optimal outcome} designated by $\lambda$.

\section{Lower and upper bounds on components of implicit Pareto optimal outcomes}
\label{bounds}
The general theory of lower and upper bounds on components of implicit Pareto optimal outcomes is given in \cite{Kaliszewski_Miroforidis_2019_genbounds}.
To calculate the bounds, one needs two finite sets (that satisfy certain properties) namely a \textit{lower shell} ($S_L \subseteq X_0$) and \textit{upper shell} ($S_U \subseteq {\mathbb{R}}^n$).

Given $\lambda$, $S_L$, and $S_U$, the theory provides formulas for calculating lower and upper bounds on $f_l(x^{P_{opt}}(\lambda))$, $l=1,\ldots,k$. That is,
\begin{equation}
\label{eq_bounds}
L_l(S_L,\lambda) \leq f_l(x^{P_{opt}}(\lambda)) \leq U_l(S_U,\lambda), \,\, l=1,\ldots,k.
\end{equation}
Formulas for lower bounds $L_l(S_L,\lambda)$ and upper bounds $U_l(S_U,\lambda)$ are shown in \cite{Kal_Mir_2022_Probing}. In that work, all those elements of the theory of lower and upper bounds that are required to understand the rest of the current work are presented in a synthetic way.

In addition, and of great relevance to the current work, the theory specifies that only elements $x \in S_U$ appropriately located with respect to the vector of lower bounds $L(S_L,\lambda):=(L_1(S_L,\lambda),\ldots,L_k(S_L,\lambda))$ can provide upper bounds $U_{\bar{l}}(\{x\},\lambda)=f_{\,\bar{l}\,}(x)$ for $f_{\,\bar{l}\,}(x^{P_{opt}}(\lambda))$ for some $\bar{l}$. 
This is specified by the following lemma defined in \cite{Kaliszewski_Miroforidis_2019_genbounds}. 

\begin{lemma}
\label{proposition3.2}
Given lower shell $S_L$ and upper shell $S_U$. Suppose $x \in S_U$ and $L_{\,\bar{l}\,}(S_L, \lambda) \leq f_{\,\bar{l}\,}(x)$ for some $\bar{l}$ and $L_l(S_L, \lambda) \geq f_l(x)$ for all \ $l=1,\ldots,k, \ l \neq \bar{l}$.
Then  $x$ provides an upper bound for  $f_{\,\bar{l}\,}(x^{P_{opt}}(\lambda))$, namely $f_{\,\bar{l}\,}(x^{P_{opt}}(\lambda)) \leq f_{\,\bar{l}\,}(x)$.
\end{lemma}

Let $\bar{S}_U \subseteq S_U$ be a set of elements fulfilling Lemma \ref{proposition3.2} for some $\bar{l} \in \{1,\ldots,k\}$. If  $\bar{S}_U \neq \emptyset$ then each $x \in \bar{S}_U$ can provide an upper bound on $f_{\,\bar{l}\,}(x^{P_{opt}}(\lambda))$, and $U_{\bar{l}}(S_U,\lambda)=\min_{x \in \bar{S}_U}U_{\bar{l}}(\{x\},\lambda)$. If $\bar{S}_U = \emptyset$ then $U_{\bar{l}}(S_U,\lambda)=\hat{y}_{\bar{l}}$. In Section \ref{numerical_experiments}, we show how to set $U_{\bar{l}}(S_U,\lambda)$ when $\hat{y}_{\bar{l}}$ is not known.


Further on,  $U(S_U,\lambda):=(U_1(S_U,\lambda),\ldots,U_k(S_U,\lambda))$ denotes the vector of upper bounds.

\subsection{The interval representation of the implicit Pareto optimal outcome}
\label{interval_representation}
Given $\lambda$, $S_L$, and $S_U$, $R(S_L, S_U,\lambda):=([L_1(S_L,\lambda),U_1(S_U,\lambda)],\ldots, [L_k(S_L,\lambda),U_k(S_U,\lambda)])$ is the interval representation of  $f(x^{P_{opt}}(\lambda))$.

\begin{figure}
\includegraphics[scale=0.4]{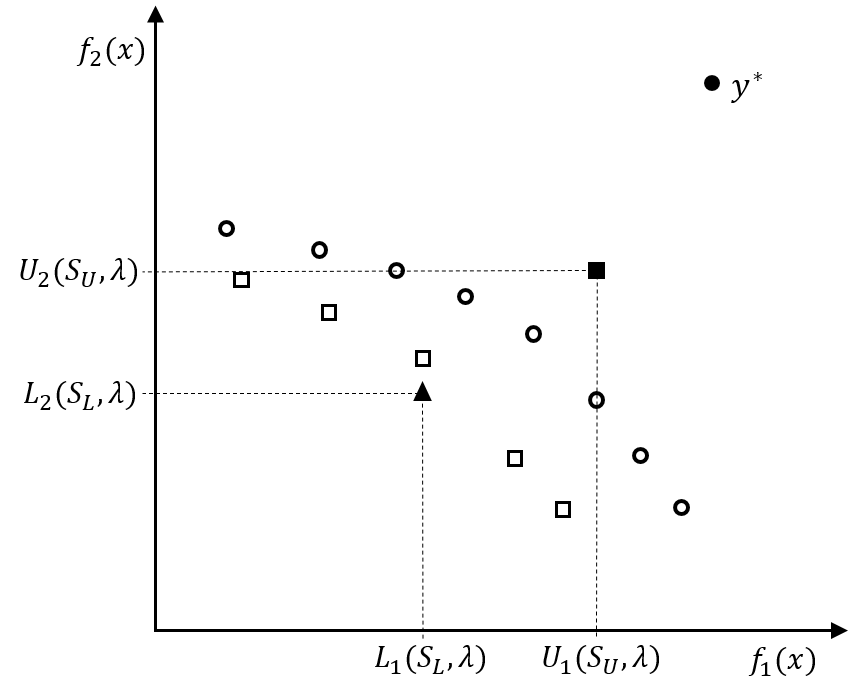} 
\centering
\caption{Components of $R(S_L, S_U,\lambda)$: $\square$, $\circ$ -- images of lower shell $S_L$ and upper shell $S_U$ elements, respectively, in the objective space, $\blacktriangle$ -- vector of lower bounds, $\blacksquare$ -- vector of upper bounds.}
\label{Fig_interval_repr}
\end{figure}

For $k=2$, components of the interval representation of $f(x^{P_{opt}}(\lambda))$, lower and upper shells as well as vectors of lower and upper bounds are illustrated in Figure \ref{Fig_interval_repr}.

To gauge the quality of $R(S_L, S_U,\lambda)$, we calculate 
\begin{equation}
\label{eq_GAP}
G_{P_{sub},l}(\lambda) := 100\times \frac{U_l(S_U, \lambda) - L_l(S_L, \lambda)}{U_l(S_U, \lambda)}, \ l=1,\ldots,k.
\end{equation}

$G_{P_{sub}}(R(S_L, S_U,\lambda)):=(G_{P_{sub},1}(\lambda),\ldots,G_{P_{sub},k}(\lambda))$ forms the \textit{Pareto suboptimality gap of interval representation} $R(S_L, S_U,\lambda)$.

\section{Providing interval representations of implicit Pareto optimal outcomes}
\label{deriv_interval_repr}

Given $\lambda$, we assume that there is a time limit $T^L$ on solving problem (\ref{eq2.2}) by a MIP solver. We also assume that if the MIP solver can not derive the solution to (\ref{eq2.2}), i.e., $x^{P_{opt}}(\lambda)$, within $T^L$, then it provides incumbent ${INC}^\lambda$ that is the approximation of $x^{P_{opt}}(\lambda)$. 
In this case, our goal is to provide $R(S_L, S_U,\lambda)$ calculated on some lower shell $S_L$ and on some upper shell $S_U$.

\subsection{The derivation of lower and upper shells}
\label{deriv_SL_SU}
As in \cite{Kal_Mir_2022_Probing}, we will use $S_L:=\{ {INC}^\lambda \}$ as a valid lower shell one can use to calculate $L_l(S_L,\lambda), l=1,\ldots,k$.

To populate upper shell $S_U$, the following two Lemmas defined in  \cite{Kal_Mir_2022_Probing} can be used.
\begin{lemma}
\label{relax_lemma}
Given $\lambda'$, solution $x'$ to the relaxation of problem (\ref{eq2.2}) with $X^{'}_0$, $X^{'}_0 \supset X_0$, is not dominated by solution $x$ to problem (\ref{eq2.2}) for any $\lambda$.
\end{lemma}

\begin{lemma}
\label{relax_lemma_MO}
If $x$ is a Pareto optimal solution to the relaxation of problem (\ref{eq2.1}) with $X'_0$, then set $\{x\}$ is an upper shell to problem (\ref{eq2.1}).
\end{lemma}
Given $X^{'}_0 \supset X_0$ and $\lambda$, let $\text{ChebRLX}(X^{'}_0, \lambda)$ denote the Chebyshev scalarization (problem (\ref{eq2.2})) of some relaxation of the MOMIP problem (problem (\ref{eq2.1})) with feasible set $X^{'}_0$ for some $\lambda$.
Based on Lemmas \ref{relax_lemma} and \ref{relax_lemma_MO}, one can derive a single-element upper shell by solving $\text{ChebRLX}(X^{'}_0, \lambda)$.
Given $X^{'}_0 \supset X_0$, the sum of such single-element upper  shells derived for different vectors $\lambda$ forms upper shell $S_U$.

The surrogate relaxation of problem (\ref{eq2.1}) with $X^{'}_0(\mu) :=\{ x \in X \, | $ $\, \sum^m_{p=1} \mu_p g_{p}(x)$ $  \leq \sum^m_{p=1} \mu_p b_p\}$ instead of $X_0$, where $\mu_p\ge 0,\,p=1,\ldots,m$, $\mu  \neq 0$, is a valid relaxation of this problem (see \cite{Kal_Mir_2022_Probing}). $\mu$ is a vector of surrogate multipliers. We will use this type of relaxation with $\mu$ as a parameter to derive elements of an upper shell.
We also follow what has been shown in \cite{Kal_Mir_2022_Probing} on large-scale instances of the  MOMIP problem that

Given $\lambda$ and $S_L$, based on Lemma 1, element $x$ of upper shell $S_U$ is a source of an upper bound on $f_{\,\bar{l}\,}(x^{P_{opt}}(\lambda))$, $\bar{l} \in \{1,\ldots,k\}$, when $f(x)$ is appropriately located with respect to the vector of lower bounds $L(S_L,\lambda)$. In \cite{Miroforidis_2021_Markowitz}, an idea of how to derive an upper shell that consists of an element useful to calculate upper bounds on $f_{\,\bar{l}\,}(x^{P_{opt}}(\lambda))$ has been proposed. This idea is to probe the objective space by perturbing components of vector $\lambda$. Yet, there is no algorithmic approach in \cite{Miroforidis_2021_Markowitz} doing that. In the current work, we try to fill in this gap.

\begin{figure}
\includegraphics[scale=0.4]{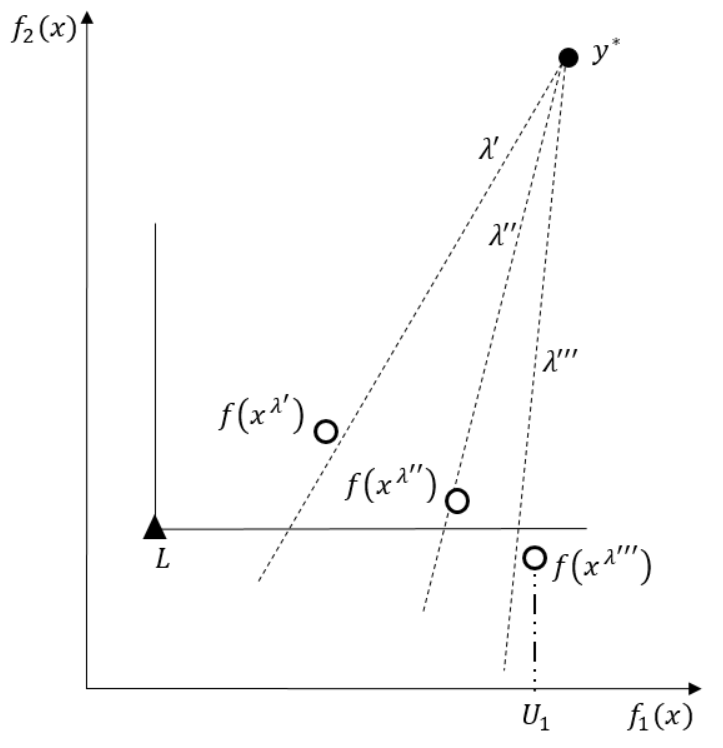} 
\centering
\caption{Deriving upper shell $S^1_U$ whose element $x^{{\lambda}^{'''}}$ is a source of an upper bound, $U_1$, for $f_{1}(x^{P_{opt}}(\lambda))$ with some $\lambda$\,: $\circ$ -- image of upper shell  $S^1_U$ in the objective space, $\blacktriangle$ -- vector of lower bounds.}
\label{Fig_FUS_Idea}
\end{figure}

For $k=2$, the idea of an algorithm for deriving upper shell $S^{1}_U$ whose some element is a source of an upper bound on $f_{1}(x^{P_{opt}}(\lambda))$ is shown in Figure \ref{Fig_FUS_Idea}.
Let us assume that vector $\mu$ is given. At the beginning, $S^{1}_U:=\emptyset$. In the first step, we set the first probing vector ${\lambda}^{'}:=({\lambda}_1+\delta,{\lambda}_2-\delta)$, ${\lambda}^{'}_2 > 0$, for some $\delta > 0$ 
(as a probing vector, we exclude $\lambda$ because we expect that the corresponding solution will not be properly located with respect to the vector of lower bounds, see \cite{Kal_Mir_2021_Cooperative}). Let $x^{{\lambda}^{'}}$ be the solution to $\text{ChebRLX}(X^{'}_{0}(\mu), {\lambda}^{'})$. $S^{1}_U:=S^{1}_U\cup\{x^{{\lambda}^{'}}\}$. Based on Lemma 1, $x^{{\lambda}^{'}}$ is not a source of an upper bound on $f_{1}(x^{P_{opt}}(\lambda))$ because $f(x^{{\lambda}^{'}})$ is not appropriately located with respect to the vector of lower bounds $L=L(S_L, \lambda)$.
Hence, we continue. In the second step, we set ${\lambda}^{''}:=({\lambda}^{'}_1+\delta,{\lambda}^{'}_2-\delta)$,  
${\lambda}^{''}_2 > 0$.
Let $x^{{\lambda}^{''}}$ be the solution to  $\text{ChebRLX}(X^{'}_{0}(\mu), {\lambda}^{''})$. 
$S^{1}_U:=S^{1}_U\cup\{x^{{\lambda}^{''}}\}$. Based on Lemma 1, $x^{{\lambda}^{''}}$ is not a source of an upper bound on $f_{1}(x^{P_{opt}}(\lambda))$.
Hence, we continue. In the third step, we set ${\lambda}^{'''}:=({\lambda}^{''}_1+\delta,{\lambda}^{''}_2-\delta)$,  
${\lambda}^{'''}_2 > 0$.
Let $x^{{\lambda}^{'''}}$ be the solution to $\text{ChebRLX}(X^{'}_{0}(\mu), {\lambda}^{'''})$. 
$S^{1}_U:=S^{1}_U\cup\{x^{{\lambda}^{'''}}\}$. Based on Lemma 1, $x^{{\lambda}^{'''}}$ is a source of an upper bound on $f_{1}(x^{P_{opt}}(\lambda))$. So, $U_1 = U(S^1_U,\lambda) = f_1(x^{{\lambda}^{'''}})$. As elements $x^{{\lambda}^{'}}$, $x^{{\lambda}^{''}}$, and $x^{{\lambda}^{'''}}$ are Pareto optimal solutions to the relaxation of the MOMIP problem with $X^{'}_{0}(\mu)$,  $S^{1}_U$ is a valid upper shell. To obtain an upper bound on $f_{2}(x^{P_{opt}}(\lambda))$, we need to derive upper shell $S^{2}_U$.
To do this, in the first step, we set ${\lambda}^{'}:=({\lambda}_1-\delta,{\lambda}_2+\delta)$,  ${\lambda}^{'}_1 > 0 $, and proceed in the same way.

Given $\bar{l} \in \{1,\ldots,k\}$, the  FindUpperShell  algorithm tries to derive upper shell $S^{\bar{l}}_U$ whose element $x$ is a source of an upper bound on $f_{\,\bar{l}\,}(x^{P_{opt}}(\lambda))$, i.e. $U_{\bar{l}}(S^{\bar{l}}_U, \lambda)$.
\begin{algorithm}
\NoCaptionOfAlgo
\scriptsize
\caption{FindUpperShell} 
\label{FindUpperShell_l} 
\textbf{INPUT}: $\bar{l}$, $\lambda$, $y^*$, $L(S_L,\lambda)$, $\gamma>0$, $\mu$ \\
\textbf{OUTPUT}: $S^{\bar{l}}_U$,  $U_{\bar{l}}(S^{\bar{l}}_U, \lambda)$\\
\vspace{0.2cm}
\nl $S^{\bar{l}}_U := \emptyset$ ;  $U_{\bar{l}}(S^{\bar{l}}_U, \lambda) := y^{*}_{\bar{l}}$ ; \textit{Comment: initialization.}\\
\nl 
$\delta := (1-{\lambda}_{\bar{l}})/\gamma$ ; \textit{Comment: set the step size.}\\
\nl
$\lambda^{'} := \lambda$ ; $\lambda^{'}_{\bar{l}} := \lambda^{'}_{\bar{l}} + \delta$ ; $f := \text{FALSE}$ ;\\
\nl \While{$\lambda^{'}_{\bar{l}} < 1$}
{
\nl \ForEach{$l=1,\ldots,k$, $l \neq \bar{l}$}
{
	\nl
	\eIf{$\lambda^{'}_{l} - \frac{\delta}{k-1} < 0$}{
		\nl	 $f := \text{TRUE}$ ; \\
		\nl	 \textbf{break} ; \textit{Comment: break foreach.}\\
	}
	{
		\nl $\lambda^{'}_{l} := \lambda^{'}_{l} - \frac{\delta}{k-1}$ ; \\
	}
}
\nl \If{\normalfont{$f = \text{TRUE}$}}{\textbf{break} ; \textit{Comment: break while.}}

\nl Let $x$ be the solution to $\text{ChebRLX}(X^{'}_{0}(\mu), {\lambda}^{'})$ ; \\
\nl $S^{\bar{l}}_U := S^{\bar{l}}_U \cup \{x\}$ ; \\
\nl \If{$x$ \normalfont{fulfills conditions of Lemma 1 for $\bar{l}$, $L(S_L,\lambda)$, and $S_U = \{x\}$}}
{	
	$U_{\bar{l}}(S^{\bar{l}}_U, \lambda) := f_{\bar{l}}(x)$ ;
	\textbf{break} ; \textit{Comment: break while.}
}
\nl
$\lambda^{'}_{\bar{l}} := \lambda^{'}_{\bar{l}} + \delta$ ; \\
}
\vspace{0.2cm}
\nl \textbf{RETURN}: $S^{\bar{l}}_U$, $U_{\bar{l}}(S^{\bar{l}}_U, \lambda)$
\end{algorithm}

In Line 2, we set step size $\delta$ that is used to modify components of consecutive probing vectors ${\lambda}^{'}$. Parameter $\gamma>0$ controls the step size, i.e., the greater the value of parameter $\gamma$, the denser the sampling of the objective space to search for the desired element of the upper shell. In the main loop (Lines 4-14), we populate upper shell $S^{\bar{l}}_U$ checking if its new element fulfills conditions of Lemma 1 to be a valid source for  the upper bound on $f_{\bar{l}}(x^{P_{opt}}(\lambda))$.
The algorithm stops when ${\lambda}^{'}_{\bar{l}} \geq 1$ \underline{OR} some component of ${\lambda}^{'}$ is negative \underline{OR} an element that fulfills conditions of Lemma 1 is found.
Lines 6 and 9 guarantee that  $\sum^k_{l=1} {{\lambda}^{'}_l} = 1$. The exit condition of the "while" loop ensures that ${\lambda}^{'}_l > 0, l=1,\ldots,k$.

If no element of $S^{\bar{l}}_U$ satisfies conditions of Lemma 1, the algorithm simply returns the upper shell as well as the only available upper bound on $f_{\,\bar{l}\,}(x^{P_{opt}}(\lambda))$, namely $ y^{*}_{\bar{l}}$.
Otherwise, an upper bound better than  $ y^{*}_{\bar{l}}$ is returned as well as the upper shell.

\subsection{Calculating interval representations}
\label{algorithm_CHUTE}

Based on the above elements, an interval representation of the Pareto optimal outcome given by vector $\lambda$ can be calculated with the use of the Chute algorithm. Along with the interval representation, this algorithm also returns lower and upper shells that were determined during its operation. In Subsection \ref{discussion}, we explain why the algorithm also returns lower and upper shells.

Line 4 of the algorithm needs clarification. Vector $\mu$ can be set as shown in \cite{Kal_Mir_2022_Probing}, namely by taking $\mu_p := 1$, $p=1,\ldots,m$. In that work, all surrogate multipliers have the same value.
We call this version of the Chute algorithm Chute1.

Yet, in \cite{Kal_Mir_2022_Probing}, Section "Final remarks", it has been suggested that \emph{"Tighter bounds might be obtained with other values of the multipliers. This possibility is worth exploring in future works."}.
Unfortunately, there is no idea there how to select a vector of surrogate multipliers other than $(1,\ldots,1) \in \mathbb{R}^{m}$.
However, we can use the theory of duality for this purpose.

Given $\mu$, $\lambda$, let $x$ be the solution to $\text{ChebRLX}(X^{'}_{0}(\mu), {\lambda})$, and $s(\mu)$ be the objective function value of $x$. 
Based on Lemmas 2 and 3, $\{x\}$ is a valid upper shell. Let $s$ be the objective function value of the solution to problem (\ref{eq2.2}) with $\lambda$ and $X_0$.
It is a well-known fact (see, e.g., \cite{Glover_1965}, \cite{Glover_1968}) that $s \geq s(\mu)$. Hence, for a given $\lambda$ and $\mu$,  $s(\mu)$ is a lower bound on values of $s$.

Given $\lambda$, the best (highest) lower bound $s^*$ on values of $s$ is the objective function value of the solution $\mu^{*}$ to the following surrogate dual problem
\begin{equation}
\label{surr_rlx_dual}
\begin{array}{ll}
\sup_{\mu \geq 0,\mu \neq 0}\{\min_{x \in X^{'}_0(\mu)} & \max_l \lambda_l (y^*_l - f_l(x)) + \rho e^k(y^* - f(x))\} 
\end{array}
\end{equation}
that is connected to the Chebyshev scalarization (problem (\ref{eq2.2})).
Solving (\ref{surr_rlx_dual}) to optimality can be time-consuming.
Yet, a suboptimal vector of multipliers $\tilde \mu$ can be determined instead of $\mu^{*}$. It can be done with the help of a quasi-subgradient-like algorithm (we shall call it Suboptimal) by Dyer (\cite{Dyer_1980}) with the following stopping condition.\\
\emph{"Number of iterations without improving the value of the objective function in problem (\ref{surr_rlx_dual}) is greater than $N$"} \underline{OR} \emph{"time limit on optimization is greater than $T^S$ seconds".}

In the current work, we set time limits on computation, hence the above stopping condition is justified in practice.

We will use vector $\frac{(1,\ldots,1)}{||(1,\ldots,1)||} \in \mathbb{R}^{m}$, where $||.||$ is the Euclidean norm, as an initial vector of surrogate multipliers in the Suboptimal algorithm.
Under the above assumptions, this algorithm has three parameters, namely $\lambda$, $N$, and $T^S$.

A version of the Chute algorithm that uses (in Line 4) the Suboptimal algorithm to set vector of surrogate multipiers $\mu$ for a given $\lambda$, we shall call Chute2. It has two additional input parameters $N$ and $T^S$.

Let us note, that in the Chute2 algorithm, we set the vector of surrogate multipliers once for a given $\lambda$. The FindUpperShell algorithm uses perturbations of the $\lambda$ vector to sample the objective space, and for all these perturbations the same vector $\mu$ is used. It is our heuristic assumption that even using the same vector $\mu$ for various vectors $\lambda^{'}$, that are close to $\lambda$, the Chute2 algorithm is able to find a better $R(S_L, S_U,\lambda)$, so tighter upper bounds on components of $f(x^{P_{opt}}(\lambda))$, than the Chute1 algorithm.
However, it is at the cost of increasing the computation time relative to Chute1 by at most $T^S$.
We will check it experimentally in the next section.

\begin{figure}
\includegraphics[scale=0.4]{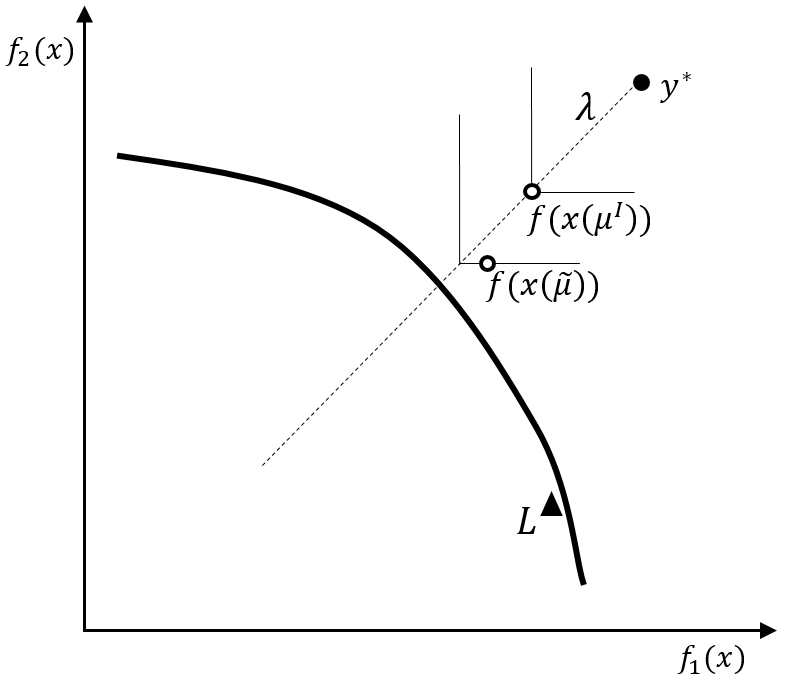} 
\centering
\caption{The idea of deriving upper shell $\{x(\tilde{\mu})\}$ whose element is a source of a better upper bound on $f_{2}(x^{P_{opt}}(\lambda))$ than the element of upper shell $\{x(\mu^{I})\}$.}
\label{Fig_Surrogate}
\end{figure}

The idea (for $k=2$ and $\rho=0$) of using suboptimal values of surrogate multipliers to get an upper shell that is a source of a better upper bound is illustrated in Figure \ref{Fig_Surrogate}. Let $\mu^I:=(1,\ldots,1) \in \mathbb{R}^{m}$, and $\tilde{\mu}$ be the output of the Suboptimal algorithm  (with some $N$ and $T^S$)  for some $\lambda$.
$x(\mu^I)$ is the solution to optimization problem (\ref{eq2.2}) with $X_0$ replaced with $X^{'}_0(\mu^I)$, and $s(\mu^{I})$ is the objective function value of $x(\mu^I)$. $x(\tilde{\mu})$ is the solution to optimization problem (\ref{eq2.2}) with $X_0$ replaced with $X^{'}_0(\tilde{\mu})$, and $s(\tilde{\mu})$ is the objective function value of $x(\tilde{\mu})$.
Vector of lower bounds $L$ is marked with a triangle. Both elements $f(x(\tilde{\mu}))$ and $f(x(\mu^{I}))$ are appropriately located with regards to $L$ (see Lemma 1) to be sources for an upper bound on $f_{2}(x^{P_{opt}}(\lambda))$.
As $s(\tilde{\mu}) > s(\mu^{I})$ (contours of the Chebyshev metric for both values  $s(\tilde{\mu})$ and $s(\mu^{I})$ are shown by solid thin lines) as well as $f_1(x(\tilde{\mu})) < f_1(x(\mu^{I}))$ \underline{AND} $f_2(x(\tilde{\mu})) < f_2(x(\mu^{I}))$, element $x(\tilde{\mu})$ is a source of a better upper bound on $f_{2}(x^{P_{opt}}(\lambda))$ for than element $x(\mu^{I})$, as upper bounds are calculated with the use of components of the upper shell elements.  In Figure \ref{Fig_Surrogate}, $f(x(\tilde{\mu}))$ is closer to the (unknown) Pareto front (represented by the solid curve) than element $f(x(\mu^{I}))$. It could happen that condition $f_1(x(\tilde{\mu})) < f_1(x(\mu^{I}))$  \underline{AND} $f_2(x(\tilde{\mu})) < f_2(x(\mu^{I}))$ does not hold. In this case, if $f_2(x(\tilde{\mu})) \geq f_2(x(\mu^{I}))$, we obtain no better upper bound on $f_{2}(x^{P_{opt}}(\lambda))$. On the other hand, if  $f_1(x(\tilde{\mu})) \geq f_1(x(\mu^{I}))$ and still $f(x(\tilde{\mu}))$ is appropriately located with regards to $L$ (see Lemma 1) to be a source for an upper bound on $f_{2}(x^{P_{opt}}(\lambda))$, we get a better upper bound on $f_{2}(x^{P_{opt}}(\lambda))$.

\begin{algorithm}
\NoCaptionOfAlgo
\scriptsize
\caption{Chute} 
\label{Chute} 
\textbf{INPUT}: $\lambda$, $y^*$, $T^L$, $\gamma>0$  \\
\textbf{OUTPUT}: $S_L$, $S_U$,  $R(S_L, S_U,\lambda)$\\
\vspace{0.2cm}
\nl Let ${INC}^{\lambda}$ be solution to problem (\ref{eq2.2}) derived within $T^L$ ;\\
\nl $S_L := \{{INC}^{\lambda}\}$ ; $S_U := \emptyset$ ;\\
\nl Let $L(S_L, \lambda)$ be the vector of lower bounds ;\\
\nl Set vector of multipliers $\mu$ ;\\
\nl \ForEach{$l=1,\ldots,k$}
{
$(S^l_U, U_l) :=\, 
$FindUpperShell($l,\lambda,y^{*},L(S_L, \lambda),\gamma, \mu)$ ; \\
$S_U := S_U \cup S^l_U$ ;\\
}
\nl $R(S_L, S_U,\lambda):=([L_1(S_L,\lambda),U_1],\ldots, [L_k(S_L,\lambda),U_k])$ ;\\
\vspace{0.2cm}
\nl \textbf{RETURN}: $S_L$, $S_U$, $R(S_L, S_U,\lambda)$ \\
\end{algorithm}

\section{Computational experiments}
\label{numerical_experiments}
In this section, we present the results of two experiments where we apply algorithms Chute1 and Chute2 presented in Subsection \ref{algorithm_CHUTE} to selected instances of the Multi-Objective Multidimensional 0-1 Knapsack Problem (MOMKP) with two and three objective functions.
The instances are demanding for modern MIP solvers.
\subsection{Multi-objective multidimensional 0-1 knapsack problem}
\label{MOMKP}

For $k > 1$, the MOMKP is formulated in the following way.
\begin{equation}
\label{eq6.1}	
\begin{array}{ll}
\text{vmax} & \left\{ 
\begin{array}{l}
	f_1(x):=\sum^n_{j=1} c_{1,j}x_j \\
	\dots \\
	f_k(x):=\sum^n_{j=1} c_{k,j}x_j
\end{array} \right. \vspace{2mm}\\
\text{s.t.} &
x \in X_0 := \{ x \, | \, \sum^n_{j=1} a_{p,j}x_j \leq b_p, \ p=1,\ldots,m, \, \vspace{1mm}\\& x_j \in \{ 0,1 \} \,, \ j=1,\ldots,n \},
\end{array}
\end{equation}
where all $a_{p,j}\,, \ c_{p,j}$ are non-negative.
In \cite{Kal_Mir_2022_Probing}, it has been explained why the MOMKP is ${\cal N}{\cal P}$-hard.

\subsection{Test instances of the MOMIP problem}
\label{test_MOMKP}
As tri-criteria instances of the MOMKP, we take two instances from \cite{Kal_Mir_2022_Probing} that were generated based on the 1st problem of the 6th group ($n=500, m=10$) of multidimensional 0-1 knapsack problems as well as on the 1st problem of the 9th group ($n=500, m=30$) (both single-objective problems are stored in Beasley OR-Library, http://people.brunel.ac.uk/$\sim$mastjjb/jeb/info.html). We call these tri-criteria instances Three6.1 and Three9.1, respectively.

By removing the third objective function of problem Three6.1, we create a bi-criteria instance called Bi6.1. Analogously, by removing the third objective function of problem Three9.1, we create a bi-criteria instance called Bi9.1.
%
%
%

Bi6.1, Bi9.1, Three6.1, and Three9.1 are our test instances of the MOMIP problem \footnote{The instances can be made available to the reader by e-mail upon request.}.

\subsection{Experimental setting}
\label{exp_setting}
Gurobi (version 10.0.0) for Windows 10 (x64) is our selected MIP solver. The optimizer is installed on the Intel Core i7-7700HQ-based laptop with 16 GB RAM.

To be consistent with limiting optimization time, we do not derive element $\hat y$ to optimality.
So, it is not known.
Instead, we separately maximize each objective function of problem (\ref{eq6.1}) within the time limit equal to 400 seconds. For instances Three6.1 and Three9.1, we set $y^*_l,\ l=1,2,3$, to the best upper bound (provided by the MIP solver) on values of the respective objective functions (none of these maxima the MIP solver determined in this time limit). Thus, for Three6.1 $y^* := (128872,131116,131738)$, and for Three9.1 $y^* := (119379.88,119365,118122)$.
Obviously, ${\hat y}_l < y^*_l$.
Hence, such $y^*$ approximating  $\hat y$ can be used in (\ref{eq2.2}) as well as when calculating lower and upper bounds.
To avoid redundant calculations, for instances Bi6.1 and Bi9.1, we take only the first two components of respective $y^{*}$.

We set $T^L := 1200$ seconds, $\rho := 0.001$.
The absolute lower bound on values of all objective functions of the MOMKP is 0 (see Subsection \ref{MOMKP}), so this value is used for calculating lower bounds $L_l(S_L,\lambda), l=1,\ldots,k.$ (for details, see \cite{Kal_Mir_2022_Probing}).

For instances Bi6.1 and Bi9.1, we generate one set of five vectors $\lambda$ uniformly sampled from two-dimensional unit simplex (see \cite{Smith_Tromble_2004}), and obtain corresponding to them vectors of lower bounds $L(S_L,\lambda)$ (Table \ref{Table_lambdas_Bi_61_91}). 
\begin{table}[htbp]
\centering
\setlength{\tabcolsep}{2pt}
\begin{tabular}[c]{c|cc|cc|cc}
\multicolumn{ 1}{c}{No} &
\multicolumn{ 2}{c}{$\lambda$} &
\multicolumn{ 2}{c}{$L(S_L,\lambda)$ for Bi6.1}&
\multicolumn{ 2}{c}{$L(S_L,\lambda)$ for Bi9.1}\\
\hline
1 & 0.055 & 0.945 & 114253.29 & 130251.56 & 104466.45 & 118482.17\\ 
2 & 0.116 & 0.884 & 116707.61 & 129508.69 & 107215.43 & 117756.82\\ 
3 & 0.733 & 0.267 & 125690.15 & 122399.79 & 116288.83 & 110899.20\\ 
4 & 0.397 & 0.603 & 122075.81 & 126638.06 & 112806.06 & 115033.24\\ 
5 & 0.439 & 0.561 & 122514.80 & 126139.05 & 113385.01 & 114671.51\\
\end{tabular}
\caption{Vectors $\lambda$ and lower bounds for test problems Bi6.1 and Bi9.1.}
\label{Table_lambdas_Bi_61_91}
\end{table}
For instances Three6.1 and Three9.1, we generate separate sets of five vectors $\lambda$, uniformly sampled from the three-dimensional unit simplex, and obtain corresponding to them vectors of lower bounds $L(S_L,\lambda)$ (Tables \ref{Table_lambdas_Tri_61} and \ref{Table_lambdas_Tri_91}, respectively) .
\begin{table}[htbp]
\centering
\setlength{\tabcolsep}{2pt}
\begin{tabular}[c]{c|ccc|ccc}
\multicolumn{ 1}{c}{No} &
\multicolumn{ 3}{c}{$\lambda$} &
\multicolumn{ 3}{c}{$L(S_L,\lambda)$}\\
\hline
1 & 0.187 & 0.770 & 0.043 & 118622.54 & 128616.78 & 87944.84\\ 
2 & 0.521 & 0.324 & 0.155 & 124156.03 & 123541.43 & 115957.65\\ 
3 & 0.067 & 0.680 & 0.253 & 90480.66 & 127282.50 & 121460.00\\ 
4 & 0.359 & 0.295 & 0.346 & 120988.48 & 121527.94 & 123559.14\\ 
5 & 0.136 & 0.078 & 0.786 & 115623.82 & 108141.30 & 129431.77\\ 
\end{tabular}
\caption{Vectors $\lambda$ and lower bounds $L(S_L,\lambda)$ for test problem Three6.1.}
\label{Table_lambdas_Tri_61}
\end{table}
\begin{table}[htbp]
\centering
\setlength{\tabcolsep}{2pt}
\begin{tabular}[c]{c|ccc|ccc}
\multicolumn{ 1}{c}{No} &
\multicolumn{ 3}{c}{$\lambda$} &
\multicolumn{ 3}{c}{$L(S_L,\lambda)$}\\
\hline
1 & 0.351 & 0.351 & 0.298 & 111876.06 & 111861.18 & 109288.07 \\ 
2 & 0.243 & 0.143 & 0.614 & 110262.24 & 103915.65 & 114504.59\\ 
3 & 0.278 & 0.494 & 0.228 & 110549.31 & 114387.77 & 107363.36\\ 
4 & 0.179 & 0.471 & 0.350 & 105139.63 & 113934.39 & 110819.31\\ 
5 & 0.407 & 0.014 & 0.579 & 112514.84 & 0.00 & 113292.80 \\ 

\end{tabular}
\caption{Vectors $\lambda$ and lower bounds $L(S_L,\lambda)$ for test problem Three9.1.}
\label{Table_lambdas_Tri_91}
\end{table}

For all problem instances, for no vector $\lambda$, the selected MIP solver derived the solution to problem (2) in the assumed $T^L=1200$ seconds. Thus, the use of the Chute algorithm to determine interval representations of Pareto optimal outcomes given by vectors $\lambda$ is justified.  

We conduct two numerical experiments. In experiment 1, we test the behavior of algorithm Chute1. In experiment 2, we test the behavior of algorithm Chute2.
In both experiments, on generated test instances, we test the behavior of the Chute algorithm for $\gamma := 10, 30, 50$. Given $\lambda$, we check the impact of parameter $\gamma$ on components of $G_{P_{sub}}(R(S_L, S_U,\lambda))$.

In the tables with results of the experiments, the meaning of the columns is as follows.
\begin{itemize}
\item $U(S_U,\lambda)$ -- components of the vector of upper bounds;
\item ${GAP}_{P_{sub}}\%$ -- components of $G_{P_{sub}}(R(S_L, S_U,\lambda))$;
\item $|S_U|$ -- the number of elements in the upper shell;
\item Time $S_U$ (s) -- time to derive the upper shell (in seconds); for Chute2, the running time of the Suboptimal algorithm is given in parentheses.
\end{itemize}


In bold, we indicate the improvement of a single component of $G_{P_{sub}}(R(S_L, S_U,\lambda))$ when changing the value of parameter $\gamma$ to a higher value, namely, we consider changes from $\gamma=10$ to $\gamma=30$ and from $\gamma=30$ to $\gamma=50$.
By underscore, we indicate the deterioration of a single component of $G_{P_{sub}}(R(S_L, S_U,\lambda))$ when changing the value of parameter $\gamma$ to a higher value, namely, we consider changes from $\gamma=10$ to $\gamma=30$ and from $\gamma=30$ to $\gamma=50$.

\subsection{Experiment 1}
\label{exp_1}
In this experiment, we check the behavior of the Chute1 algorithm.

The results for instance Bi6.1 are shown in Tables \ref{Table_exp1_6_1_UB}--\ref{Table_exp1_6_1_Su_Time}, and the results for instance Bi9.1 are shown in Tables \ref{Table_exp1_9_1_UB}--\ref{Table_exp1_9_1__Su_Time}.
The results for instance Three6.1 are shown in Tables \ref{Table_exp1_3crit_6_1_UB}--\ref{Table_exp1_3crit_6_1_Su_Time}, and the results for instance Three9.1 are shown in Tables \ref{Table_exp1_3crit_9_1_UB}--\ref{Table_exp1_3crit_9_1_Su_Time}.

For instance Bi6.1, when changing $\gamma=10$ to $\gamma=30$, we observe an improvement in at least one component of $G_{P_{sub}}(R(S_L, S_U,\lambda))$ for four vectors $\lambda$, although only in two cases ($\lambda$ no. 3 and 5) two components improve. Yet, when changing $\gamma=30$ to $\gamma=50$, we observe an improvement  of $G_{P_{sub}}(R(S_L, S_U,\lambda))$ in at least one of its components for three vectors $\lambda$.
For $\lambda$ no. 3, we observe a deterioration of the first component of $G_{P_{sub}}(R(S_L, S_U,\lambda))$, and at the same time, an improvement of the second.
For $\lambda$ no. 5, we observe an improvement of the first component of $G_{P_{sub}}(R(S_L, S_U,\lambda))$, and at the same time, a deterioration of the second.

When changing $\gamma=10$ to $\gamma=50$, we observe an improvement in at least one component of $G_{P_{sub}}(R(S_L, S_U,\lambda))$ for all vectors $\lambda$.

For instance Bi9.1, we observe a similar phenomenon (an improvement as well as a deterioration), although when changing $\gamma=10$ to $\gamma=30$, we observe an improvement in only one component of $G_{P_{sub}}(R(S_L, S_U,\lambda))$ for just two vectors $\lambda$. When changing $\gamma=10$ to $\gamma=50$, we observe an improvement in $G_{P_{sub}}(R(S_L, S_U,\lambda))$ for vectors $\lambda$ no. 1--4 (at least one component improves). For $\lambda$ no. 5, $G_{P_{sub}}(R(S_L, S_U,\lambda))$ remains unchanged.

For instances Bi6.1 and Bi9.1, the higher the value of parameter $\gamma$ (higher sampling density of the objective space 
), the more numerous  the derived upper shells are.  For each vector $\lambda$, the time to derive the corresponding upper shell is a small fraction of the assumed time limit $T^L=1200$ seconds.

Let us check the results for tri-criteria instances. For instance Three6.1, when changing $\gamma=10$ to $\gamma=30$, we observe an improvement 
of $G_{P_{sub}}(R(S_L, S_U,\lambda))$ in at least one of its components for three vectors $\lambda$, although only for one ($\lambda$ no. 4) all components improve. When changing $\gamma=30$ to $\gamma=50$, we observe an improvement  of $G_{P_{sub}}(R(S_L, S_U,\lambda))$ in at least one of its components for three vectors $\lambda$. 
We also observe a deterioration of the first component of $G_{P_{sub}}(R(S_L, S_U,\lambda))$ for $\lambda$ no. 2, and, at the same time, an improvement on its third component.
When changing $\gamma=10$ to $\gamma=50$, we observe an improvement in $G_{P_{sub}}(R(S_L, S_U,\lambda))$ for all vectors $\lambda$ (at least one component improves). For instance Three9.1 and vectors $\lambda$ no. 2--5, upper bounds on all components of $f(x^{P_{opt}}(\lambda))$ are equal to the corresponding component of $y^{*}$.  Only for $\lambda$ no. 1, $f_3(x^{P_{opt}}(\lambda))<y_3^{*}$, and when changing $\gamma=10$ to $\gamma=30$ as well as $\gamma=30$ to $\gamma=50$, there is an improvement  only for the third component of $G_{P_{sub}}(R(S_L, S_U,\lambda))$.
For instances Three6.1 and Three9.1, the higher the value of parameter $\gamma$ (higher sampling density of the objective space 
), the more numerous the derived upper shells are. For instance Three6.1 and instance Three9.1 with $\gamma=10$, for most vectors  $\lambda$, the time to derive the corresponding upper shell is a small fraction of the assumed time limit $T^L=1200$ seconds. However, for instance Three9.1 with $\gamma=30$ and $\gamma=50$, the time to derive the corresponding upper shell increases significantly compared to $\gamma=10$, for all vectors $\lambda$. 

We checked that for all instances, for no vector $\lambda$, the MIP solver derived the optimal solution to problem (\ref{eq2.2}) within time limit $T^L+ \text{Time}\, S_U$.


\begin{table}[htbp]
\centering
\setlength{\tabcolsep}{2pt}
\begin{tabular}[c]{c|cc|cc|cc}
\multicolumn{ 1}{c}{} & 
\multicolumn{ 6}{c}{$U(S_U,\lambda)$} \\
\multicolumn{ 1}{c}{No} & 
\multicolumn{ 2}{c}{$\gamma=10$} & 
\multicolumn{ 2}{c}{$\gamma=30$} & 
\multicolumn{ 2}{c}{$\gamma=50$} \\ 
\hline
1 & 120964.00 & 131117.00 & 120466.00 & 131117.00 & 120093.00 & 131117.00 \\
2 & 121666.00 & 131117.00 & 121666.00 & 131117.00 & 121441.00 & 131117.00\\ 
3 & 127790.00 & 126078.00 & 127635.00 & 125842.00 & 127646.00 & 125642.00\\ 
4 & 125252.00 & 128990.00 & 124924.00 & 128990.00 & 124852.00 & 128990.00\\ 
5 & 125502.00 & 128779.00 & 125335.00 & 128665.00 & 125307.00 & 128710.00\\ 
\end{tabular}
\caption{Chute1, vectors of upper bounds for test problem Bi6.1 and $\gamma \in{\{10,30,50\}}$}
\label{Table_exp1_6_1_UB}
\end{table}

\begin{table}[htbp]
\centering
\setlength{\tabcolsep}{2pt}
\begin{tabular}[c]{c|cc|cc|cc}
\multicolumn{ 1}{c}{} & 
\multicolumn{ 6}{c}{$GAP_{P_{sub}}\%$} \\
\multicolumn{ 1}{c}{No} & 
\multicolumn{ 2}{c}{$\gamma=10$} & 
\multicolumn{ 2}{c}{$\gamma=30$} & 
\multicolumn{ 2}{c}{$\gamma=50$} \\ 
\hline
1 & 5.55 & 0.66  & \textbf{5.16} & 0.66 & \textbf{4.86} & 0.66 \\
2 & 4.08 & 1.23  & 4.08 & 1.23 & \textbf{3.90} & 1.23 \\  
3 & 1.64 & 2.92 & \textbf{1.52} & \textbf{2.74} & \underline{1.53} & \textbf{2.58} \\  
4 & 2.54 & 1.82 & \textbf{2.28} & 1.82 & \textbf{2.22} & 1.82 \\  
5 & 2.38 & 2.05 & \textbf{2.25} & \textbf{1.96} & \textbf{2.23} & \underline{2.00}\\  
\end{tabular}
\caption{Chute1, $GAP_{P_{sub}}\%$ for test problem Bi6.1 and $\gamma \in{\{10,30,50\}}$}
\label{Table_exp1_6_1_GAP}
\end{table}

\begin{table}[htbp]
\centering
\setlength{\tabcolsep}{2pt}
\begin{tabular}[c]{c|c|c|c|c|c|c}
\multicolumn{ 1}{c}{} & 
\multicolumn{ 3}{c}{$|S_U|$} & 
\multicolumn{ 3}{c}{Time$\,S_U (s)$}\\
\multicolumn{ 1}{c}{No} & 
\multicolumn{ 1}{c}{$\gamma=10$} & 
\multicolumn{ 1}{c}{$\gamma=30$} & 
\multicolumn{ 1}{c}{$\gamma=50$} & 
\multicolumn{ 1}{c}{$\gamma=10$} & 
\multicolumn{ 1}{c}{$\gamma=30$} & 
\multicolumn{ 1}{c}{$\gamma=50$}  \\ 
\hline
1 & 11  & 32 & 52 & 1.90 & 8.00 & 11.49\\
2 & 11 & 33 & 54 & 2.02 & 5.97 & 11.52\\  
3 & 8 & 21 & 34 & 2.07 & 4.26 & 9.12\\  
4 & 7 & 19 & 31 &2.38 & 5.90 & 9.50\\  
5 & 7 & 19 & 32 & 2.08 & 5.63 & 8.59\\  

\end{tabular}
\caption{Chute1, values of $|S_U|$, and Time$\,S_U (s)$ for test problem Bi6.1 and $\gamma \in{\{10,30,50\}}$}
\label{Table_exp1_6_1_Su_Time}
\end{table}

\begin{table}[htbp]
\centering
\setlength{\tabcolsep}{2pt}
\begin{tabular}[c]{c|cc|cc|cc}
\multicolumn{ 1}{c}{} & 
\multicolumn{ 6}{c}{$U(S_U,\lambda)$} \\
\multicolumn{ 1}{c}{No} & 
\multicolumn{ 2}{c}{$\gamma=10$} & 
\multicolumn{ 2}{c}{$\gamma=30$} & 
\multicolumn{ 2}{c}{$\gamma=50$} \\ 
\hline
1 & 116289.00 & 119365.00 & 116289.00 & 119365.00 & 116193.00 & 119365.00 \\
2 & 117277.00 & 119365.00 & 117050.00 & 119365.00 & 117057.00 & 119365.00\\ 
3 & 119379.88 & 118456.00 & 119379.88 & 118456.00 &  119379.88 & 118404.00 \\ 
4 & 119379.88 & 119365.00 & 119295.00 & 119365.00 & 119329.00 & 119365.00 \\ 
5 & 119379.88 & 119365.00 & 119379.88 & 119365.00 & 119379.88 & 119365.00 \\
\end{tabular}
\caption{Chute1, vectors of upper bounds for test problem Bi9.1 and $\gamma \in{\{10,30,50\}}$}
\label{Table_exp1_9_1_UB}
\end{table}

\begin{table}[htbp]
\centering
\setlength{\tabcolsep}{2pt}
\begin{tabular}[c]{c|cc|cc|cc}
\multicolumn{ 1}{c}{} & 
\multicolumn{ 6}{c}{$GAP_{P_{sub}}\%$} \\ 
\multicolumn{ 1}{c}{No} & 
\multicolumn{ 2}{c}{$\gamma=10$} & 
\multicolumn{ 2}{c}{$\gamma=30$} & 
\multicolumn{ 2}{c}{$\gamma=50$} \\ 
\hline
1 & 10.17 & 0.74 & 10.17 & 0.74 & \textbf{10.09} & 0.74 \\ 
2 & 8.58 & 1.35 & \textbf{8.40} & 1.35 & \underline{8.41} & 1.35 \\   
3 & 2.59 & 6.38 & 2.59 & 6.38 & 2.59 & \textbf{6.34} \\    
4 & 5.51 & 3.63 & \textbf{5.44} & 3.63 & \underline{5.47} & \underline{3.63} \\  
5 & 5.02 & 3.93 & 5.02 & 3.93 & 5.02 & 3.93 \\  

\end{tabular}
\caption{Chute1, $GAP_{P_{sub}}\%$ for test problem Bi9.1 and $\gamma \in{\{10,30,50\}}$}
\label{Table_exp1_9_1_GAP_Su_Time}
\end{table}

\begin{table}[htbp]
\centering
\setlength{\tabcolsep}{2pt}
\begin{tabular}[c]{c|c|c|c|c|c|c}
\multicolumn{ 1}{c}{} & 
\multicolumn{ 3}{c}{$|S_U|$} & 
\multicolumn{ 3}{c}{Time$\,S_U (s)$}\\
\multicolumn{ 1}{c}{No} & 
\multicolumn{ 1}{c}{$\gamma=10$} & 
\multicolumn{ 1}{c}{$\gamma=30$} & 
\multicolumn{ 1}{c}{$\gamma=50$} & 
\multicolumn{ 1}{c}{$\gamma=10$} & 
\multicolumn{ 1}{c}{$\gamma=30$} & 
\multicolumn{ 1}{c}{$\gamma=50$}  \\ 
\hline
1 & 12 & 36 & 59 & 3.06 & 11.22 & 20.44\\ 
2 & 14 & 40 & 67 & 3.91 & 10.70 & 17.21\\   
3 & 17 & 51 & 84 & 3.82 & 13.05 & 21.78\\    
4 & 20 & 4.77 & 99 & 59 & 13.79 & 24.32\\  
5 & 20 & 4.54 & 100& 60 & 14.84 & 21.51\\  

\end{tabular}
\caption{Chute1, values of $|S_U|$, and Time$\,S_U (s)$ for test problem Bi9.1 and $\gamma \in{\{10,30,50\}}$}
\label{Table_exp1_9_1__Su_Time}
\end{table}

\begin{table}[htbp]
\centering
\setlength{\tabcolsep}{2pt}
\begin{tabular}[c]{c|ccc|ccc|ccc}
\multicolumn{ 1}{c}{} & 
\multicolumn{ 9}{c}{$U(S_U,\lambda)$} \\
\multicolumn{ 1}{c}{No} & 
\multicolumn{ 3}{c}{$\gamma=10$} & 
\multicolumn{ 3}{c}{$\gamma=30$} & 
\multicolumn{ 3}{c}{$\gamma=50$} \\ 
\hline
1 & 128872.00 & 131116.00 & 122914.00 & 128872.00 & 131116.00 & 122471.00 & 128872.00 & 131116.00 & 122387.00 \\
2 & 127143.00 & 127325.00 & 124450.00 & 127069.00 & 127325.00 & 123193.00 & 127105.00 & 127325.00 & 122899.00 \\ 
3 & 124359.00 & 131116.00 & 131738.00 & 124359.00 & 131116.00 & 131738.00 & 124217.00 & 131116.00 & 131738.00 \\ 
4 & 125062.00 & 125714.00 & 126639.00 & 124809.00 & 125327.00 & 126342.00 & 124747.00 & 125273.00 & 126099.00 \\ 
5 & 128872.00 & 120016.00 & 131738.00 & 128872.00 & 120016.00 & 131738.00 & 128872.00 & 119610.00 & 131738.00 \\
\end{tabular}
\caption{Chute1, vectors of upper bounds for test problem Three6.1 and $\gamma \in{\{10,30,50\}}$}
\label{Table_exp1_3crit_6_1_UB}
\end{table}

\begin{table}[htbp]
\centering
\setlength{\tabcolsep}{2pt}
\begin{tabular}[c]{c|ccc|ccc|ccc}
\multicolumn{ 1}{c}{} & 
\multicolumn{ 9}{c}{$GAP_{P_{sub}}\%$} \\ 
\multicolumn{ 1}{c}{No} & 
\multicolumn{ 3}{c}{$\gamma=10$} & 
\multicolumn{ 3}{c}{$\gamma=30$} & 
\multicolumn{ 3}{c}{$\gamma=50$} \\ 
\hline
1 & 7.95 & 1.91 & 28.45 & 7.95 & 1.91 & \textbf{28.19} & 7.95 & 1.91 & \textbf{28.14} \\ 
2 & 2.35 & 2.97 & 6.82 & \textbf{2.29} & 2.97 & \textbf{5.87} & \underline{2.32} & 2.97 & \textbf{5.65} \\   
3 & 27.24 & 2.92 & 7.80 & 27.24 & 2.92 & 7.80 & \textbf{27.16} & 2.92 & 7.80 \\      
4 & 3.26 & 3.33 & 2.43 & \textbf{3.06} & \textbf{3.03} & \textbf{2.20} & \textbf{3.01} & \textbf{2.99} & \textbf{2.01}\\   
5 & 10.28 & 9.89 & 1.75 & 10.28 & 9.89 & 1.75 & 10.28 & 9.59 & 1.75 \\ 
\end{tabular}
\caption{Chute1, $GAP_{P_{sub}}\%$ for test problem Three6.1 and $\gamma \in{\{10,30,50\}}$}
\label{Table_exp1_3crit_6_1_GAP}
\end{table}

\begin{table}[htbp]
\centering
\setlength{\tabcolsep}{2pt}
\begin{tabular}[c]{c|c|c|c|c|c|c}
\multicolumn{ 1}{c}{} & 
\multicolumn{ 3}{c}{$|S_U|$} & 
\multicolumn{ 3}{c}{Time$\,S_U (s)$}\\
\multicolumn{ 1}{c}{No} & 
\multicolumn{ 1}{c}{$\gamma=10$} & 
\multicolumn{ 1}{c}{$\gamma=30$} & 
\multicolumn{ 1}{c}{$\gamma=50$} & 
\multicolumn{ 1}{c}{$\gamma=10$} & 
\multicolumn{ 1}{c}{$\gamma=30$} & 
\multicolumn{ 1}{c}{$\gamma=50$}  \\ 
\hline
1 & 8 & 14 & 33& 2.69 & 4.46 & 10.48\\   
2 & 11 & 21 & 50 & 3.39 & 7.80 & 23.98\\   
3 & 11 & 21 & 49 & 5.37 & 11.84 & 22.38\\      
4 & 7 & 12 & 28 & 5.46 & 9.28 & 24.29\\    
5 & 11 & 21 & 51 & 3.27 & 6.54 & 16.37\\ 
\end{tabular}
\caption{Chute1, values $|S_U|$, and Time$\,S_U (s)$ for test problem Three6.1 and $\gamma \in{\{10,30,50\}}$}
\label{Table_exp1_3crit_6_1_Su_Time}
\end{table}

\begin{table}[htbp]
\centering
\setlength{\tabcolsep}{2pt}
\begin{tabular}[c]{c|ccc|ccc|ccc}
\multicolumn{ 1}{c}{} & 
\multicolumn{ 9}{c}{$U(S_U,\lambda)$} \\
\multicolumn{ 1}{c}{No} & 
\multicolumn{ 3}{c}{$\gamma=10$} & 
\multicolumn{ 3}{c}{$\gamma=30$} & 
\multicolumn{ 3}{c}{$\gamma=50$} \\ 
\hline
1 & 119379.88 & 119365.00 & 117516.00 & 119379.8835 & 119365.00 & 117398.00 & 119379.8835 & 119365.00 & 117362.00 \\
2 & 119379.88 & 119365.00 & 118122.00 & 119379.8835 & 119365.00 & 118122.00 & 119379.8835 & 119365.00 & 118122.00 \\ 
3 & 119379.88 & 119365.00 & 118122.00 & 119379.8835 & 119365.00 & 118122.00 & 119379.8835 & 119365.00 & 118122.00 \\ 
4 & 119379.88 & 119365.00 & 118122.00 & 119379.8835 & 119365.00 & 118122.00 & 119379.8835 & 119365.00 & 118122.00 \\ 
5 & 119379.88 & 119365.00 & 118122.00 & 119379.8835 & 119365.00 & 118122.00 & 119379.8835 & 119365.00 & 118122.00 \\

\end{tabular}
\caption{Chute1, vectors of upper bounds for test problem Three9.1 and $\gamma \in{\{10,30,50\}}$}
\label{Table_exp1_3crit_9_1_UB}
\end{table}

\begin{table}[htbp]
\centering
\setlength{\tabcolsep}{2pt}
\begin{tabular}[c]{c|ccc|ccc|ccc}
\multicolumn{ 1}{c}{} & 
\multicolumn{ 9}{c}{$GAP_{P_{sub}}\%$} \\ 
\multicolumn{ 1}{c}{No} & 
\multicolumn{ 3}{c}{$\gamma=10$} & 
\multicolumn{ 3}{c}{$\gamma=30$} & 
\multicolumn{ 3}{c}{$\gamma=50$} \\ 
\hline
1 & 6.29 & 6.29 & 7.00 & 6.29 & 6.29 & \textbf{6.91} & 6.29 & 6.29 & \textbf{6.88} \\ 
2 & 7.64 & 12.94 & 3.06 & 7.64 & 12.94 & 3.06 & 7.64 & 12.94 & 3.06 \\   
3 & 7.40 & 4.17 & 9.11 & 7.40 & 4.17 & 9.11 & 7.40 & 4.17 & 9.11\\      
4 & 11.93 & 4.55 & 6.18 & 11.93 & 4.55 & 6.18 & 11.93 & 4.55 & 6.18 \\   
5 & 5.75 & 100.00 & 4.09 & 5.75 & 100.00 & 4.09 & 5.75 & 100.00 & 4.09\\ 

\end{tabular}
\caption{Chute1, $GAP_{P_{sub}}\%$ for test problem Three9.1 and $\gamma \in{\{10,30,50\}}$}
\label{Table_exp1_3crit_9_1_GAP}
\end{table}

\begin{table}[htbp]
\centering
\setlength{\tabcolsep}{2pt}
\begin{tabular}[c]{c|c|c|c|c|c|c}
\multicolumn{ 1}{c}{} & 
\multicolumn{ 3}{c}{$|S_U|$} & 
\multicolumn{ 3}{c}{Time$\,S_U (s)$}\\
\multicolumn{ 1}{c}{No} & 
\multicolumn{ 1}{c}{$\gamma=10$} & 
\multicolumn{ 1}{c}{$\gamma=30$} & 
\multicolumn{ 1}{c}{$\gamma=50$} & 
\multicolumn{ 1}{c}{$\gamma=10$} & 
\multicolumn{ 1}{c}{$\gamma=30$} & 
\multicolumn{ 1}{c}{$\gamma=50$}  \\ 
\hline
1 & 28 & 79 & 130 & 50.09 & 163.04 & 237.71\\ 
2 & 18 & 53 & 86 & 46.24 & 324.61 & 413.15\\ 
3 & 25 & 69 & 115 & 79.80 & 222.40 & 387.76\\ 
4 & 22 & 64 & 105 & 36.09 & 104.40 & 175.58\\  
5 & 11 & 29 & 49 & 11.38 & 50.33 & 101.04\\ 

\end{tabular}
\caption{Chute1, values of $|S_U|$, and Time$\,S_U (s)$ for test problem Three9.1  and $\gamma \in{\{10,30,50\}}$}
\label{Table_exp1_3crit_9_1_Su_Time}
\end{table}

\subsection{Experiment 2}
\label{exp_2}
In this experiment, we check the behavior of the Chute2 algorithm  with parameters $N=20$ and $T^S=400$. Recall, they are related to the Suboptimal algorithm.

The results for instance Bi6.1 are shown in Tables \ref{Table_exp2_6_1_UB}--\ref{Table_exp2_6_1_Su_Time}, and the results for instance Bi9.1 are shown in Tables \ref{Table_exp2_9_1_UB}--\ref{Table_exp2_9_1_Su_Time}.
The results for instance Three6.1 are shown in Tables \ref{Table_exp2_3crit_6_1_UB}--\ref{Table_exp2_3crit_6_1_Su_Time}, and the results for instance Three9.1 are shown in Tables \ref{Table_exp2_3crit_9_1_UB}--\ref{Table_exp2_3crit_9_1_Su_Time}.

For instance Bi6.1, when changing $\gamma=10$ to $\gamma=30$, we observe an improvement in at least one component of $G_{P_{sub}}(R(S_L, S_U,\lambda))$ for all vectors $\lambda$, although only in three cases ($\lambda$ no. 3--5) two components improve. Yet, when changing $\gamma=30$ to $\gamma=50$, we observe an improvement 
of $G_{P_{sub}}(R(S_L, S_U,\lambda))$ in at least one of its components for three vectors $\lambda$ (no. 1, 2, and 5).
For $\lambda$ no. 3, we observe a deterioration of the second component of $G_{P_{sub}}(R(S_L, S_U,\lambda))$, and, at the same time, an improvement on the first one. All components of $G_{P_{sub}}(R(S_L, S_U,\lambda))$ deteriorate for $\lambda$ no. 4.
When changing $\gamma=10$ to $\gamma=50$, we observe an improvement in $G_{P_{sub}}(R(S_L, S_U,\lambda))$ for all vectors $\lambda$ (at least one component improves).

For instance Bi9.1, when changing $\gamma=10$ to $\gamma=30$, we observe an improvement in at least one component of $G_{P_{sub}}(R(S_L, S_U,\lambda))$ for all vectors $\lambda$, and for  $\lambda$ no. 2--5 both components of $G_{P_{sub}}(R(S_L, S_U,\lambda))$ improve. When changing $\gamma=30$ to $\gamma=50$, we observe an improvement in both components of $G_{P_{sub}}(R(S_L, S_U,\lambda))$ for all vectors $\lambda$ but the first one , where the first component of $G_{P_{sub}}(R(S_L, S_U,\lambda))$ deteriorates.
When changing $\gamma=10$ to $\gamma=50$, we observe an improvement in $G_{P_{sub}}(R(S_L, S_U,\lambda))$ for all vectors $\lambda$ (at least one component improves).

For instances Bi6.1 and Bi9.1, the higher the value of parameter $\gamma$ (higher sampling density of the objective space 
), the more numerous the derived upper shells are.  For instance Bi6.1, average times over all vectors $\lambda$ to derive the upper shell for $\gamma=10$, $\gamma=30$, and $\gamma=50$ are, respectively, 77.74, 85.03, and 83.28 seconds. These times are small fractions of the assumed time limit $T^L=1200$ seconds.
For instance Bi9.1, average times over all vectors $\lambda$ to derive the upper shell for $\gamma=10$, $\gamma=30$, and $\gamma=50$ are, respectively, 444.71, 535.92, and 618.98 seconds, and they are\textbf{ not small} fractions of $T^L=1200$ seconds.

Let us check the results for tri-criteria instances. For instance Three6.1, when changing $\gamma=10$ to $\gamma=30$, we observe an improvement  of $G_{P_{sub}}(R(S_L, S_U,\lambda))$ in at least one of its components for vectors $\lambda$ no. 1, and 3--5. For $\lambda$ no. 2, we observe an improvement in the third component of $G_{P_{sub}}(R(S_L, S_U,\lambda))$ as well as a deterioration of the first one. When changing $\gamma=30$ to $\gamma=50$, we observe an improvement of $G_{P_{sub}}(R(S_L, S_U,\lambda))$ in at least one of its components for vectors $\lambda$ no. 1--3, and 5. For $\lambda$ no. 4, we observe an improvement in the second and third components of $G_{P_{sub}}(R(S_L, S_U,\lambda))$ as well as a deterioration of the first one. 
When changing $\gamma=10$ to $\gamma=50$, we observe an improvement in the $G_{P_{sub}}(R(S_L, S_U,\lambda))$ for all vectors $\lambda$ (at least one component improves).

For instance Three9.1, when changing $\gamma=10$ to $\gamma=30$, we observe an improvement in all components of $G_{P_{sub}}(R(S_L, S_U,\lambda))$ for vectors $\lambda$ no. 1--4, and for $\lambda$ no. 5,
$G_{P_{sub}}(R(S_L, S_U,\lambda))$ remains unchanged. When changing $\gamma=30$ to $\gamma=50$, we observe an improvement in all components of $G_{P_{sub}}(R(S_L, S_U,\lambda))$ for vector $\lambda$ no. 1. For $\lambda$ no. 2, we observe a deterioration of all components of $G_{P_{sub}}(R(S_L, S_U,\lambda))$, and for $\lambda$ no. 3--4, we observe an improvement in two components of $G_{P_{sub}}(R(S_L, S_U,\lambda))$, and a deterioration of one of its components.
When changing $\gamma=10$ to $\gamma=50$, we observe an improvement in the $G_{P_{sub}}(R(S_L, S_U,\lambda))$ for all vectors $\lambda$ (at least one component improves) but the last one. 

For instances Three6.1 and Three9.1, the higher the value of parameter $\gamma$ (higher sampling density of the objective space 
For instance Three9.1, average times over all vectors $\lambda$ to derive the upper shell for $\gamma=10$, $\gamma=30$, and $\gamma=50$ are, respectively, 425.48, 470.58, and 531.72 seconds, and they are \textbf{not small} fractions of $T^L=1200$ seconds.

We checked that for all instances, for no vector $\lambda$, the MIP solver derived the optimal solution to problem (\ref{eq2.2}) within time limit $T^L+ \text{Time}\, S_U$.

\begin{table}[htbp]
\centering
\setlength{\tabcolsep}{2pt}
\begin{tabular}[c]{c|cc|cc|cc}
\multicolumn{ 1}{c}{} & 
\multicolumn{ 6}{c}{$U(S_U,\lambda)$} \\
\multicolumn{ 1}{c}{No} & 
\multicolumn{ 2}{c}{$\gamma=10$} & 
\multicolumn{ 2}{c}{$\gamma=30$} & 
\multicolumn{ 2}{c}{$\gamma=50$} \\ 
\hline
1 & 118423.00 & 130703.00 & 116671.00 & 130703.00 & 116230.00 & 130703.00 \\
2 & 119507.00 & 130021.00 & 118226.00 & 129980.00 & 117969.00 & 129946.00\\ 
3 & 126391.00 & 123927.00 & 126213.00 & 123235.00 & 126179.00 & 123296.00 \\ 
4 & 123079.00 & 127228.00 & 122598.00 & 127102.00 & 122657.00 & 127146.00 \\ 
5 & 123474.00 & 126864.00 & 123273.00 & 126864.00 & 123233.00 & 126766.00 \\

\end{tabular}
\caption{Chute2, upper bounds for test problem Bi6.1 and $\gamma \in{\{10,30,50\}}$}
\label{Table_exp2_6_1_UB}
\end{table}

\begin{table}[htbp]
\centering
\setlength{\tabcolsep}{2pt}
\begin{tabular}[c]{c|cc|cc|cc}
\multicolumn{ 1}{c}{} & 
\multicolumn{ 6}{c}{$GAP_{P_{sub}}\%$} \\ 
\multicolumn{ 1}{c}{No} & 
\multicolumn{ 2}{c}{$\gamma=10$} & 
\multicolumn{ 2}{c}{$\gamma=30$} & 
\multicolumn{ 2}{c}{$\gamma=50$} \\ 
\hline
1 & 3.52 & 0.35 & \textbf{2.07} & 0.35 & \textbf{1.70} & 0.35 \\
2 & 2.34 & 0.39 & \textbf{1.28} & \textbf{0.36}  & \textbf{1.07} & \textbf{0.34}\\ 
3 & 0.55 & 1.23 & \textbf{0.41} & \textbf{0.68}  & \textbf{0.39} & \underline{0.73} \\ 
4 & 0.82 & 0.46 & \textbf{0.43} & \textbf{0.37}  & \underline{0.47} & \underline{0.40} \\ 
5 & 0.78 & 0.57 & \textbf{0.62} & 0.57  & \textbf{0.58} & \textbf{0.49} \\

\end{tabular}
\caption{Chute2, $GAP_{P_{sub}}\%$ for test problem Bi6.1 and $\gamma \in{\{10,30,50\}}$}
\label{Table_exp2_6_1_GAP}
\end{table}

\begin{table}[htbp]
\centering
\setlength{\tabcolsep}{2pt}
\begin{tabular}[c]{c|c|c|c|c|c|c}
\multicolumn{ 1}{c}{} & 
\multicolumn{ 3}{c}{$|S_U|$} & 
\multicolumn{ 3}{c}{Time$\,S_U (s)$}\\
\multicolumn{ 1}{c}{No} & 
\multicolumn{ 1}{c}{$\gamma=10$} & 
\multicolumn{ 1}{c}{$\gamma=30$} & 
\multicolumn{ 1}{c}{$\gamma=50$} & 
\multicolumn{ 1}{c}{$\gamma=10$} & 
\multicolumn{ 1}{c}{$\gamma=30$} & 
\multicolumn{ 1}{c}{$\gamma=50$}  \\ 
\hline
1 & 6 & 16 & 26 & 81.16 (78.64)  & 84.79 (77.97)  & 91.44 (79.81)\\  
2 & 4 & 9 & 13 & 182.12 (180.58) & 215.77 (211.45) & 197.74 (192.33)\\ 
3 & 3 & 5 & 8 & 37.59 (36.74)   & 41.65 (38.50) & 42.82 (37.99)\\ 
4 & 2 & 3 & 6 & 42.39 (41.04)   & 41.65 (40.42) & 44.57 (41.78)\\ 
5 & 2 & 5 & 7 & 44.07 (43.33)   & 41.28 (39.66) & 39.82 (37.83)\\

\end{tabular}
\caption{Chute2, $|S_U|$, and Time$\,S_U (s)$ for test problem Bi6.1 and $\gamma \in{\{10,30,50\}}$}
\label{Table_exp2_6_1_Su_Time}
\end{table}

\begin{table}[htbp]
\centering
\setlength{\tabcolsep}{2pt}
\begin{tabular}[c]{c|cc|cc|cc}
\multicolumn{ 1}{c}{} & 
\multicolumn{ 6}{c}{$U(S_U,\lambda)$} \\
\multicolumn{ 1}{c}{No} & 
\multicolumn{ 2}{c}{$\gamma=10$} & 
\multicolumn{ 2}{c}{$\gamma=30$} & 
\multicolumn{ 2}{c}{$\gamma=50$} \\ 
\hline
1 & 110646.00 & 119365.00 & 109109.00 & 119365.00 & 109313.00 & 119365.00 \\
2 & 111592.00 & 119218.00 & 111080.00 & 119124.00 & 110985.00 & 119105.00 \\ 
3 & 118331.00 & 114263.00 & 118172.00 & 113760.00 & 118138.00 & 113655.00\\ 
4 & 115241.00 & 116795.00 & 115230.00 & 116781.00 & 115121.00 & 116732.00\\ 
5 & 115443.00 & 116518.00 & 115426.00 & 116271.00 & 115321.00 & 116233.00\\

\end{tabular}
\caption{Chute2, upper bounds for test problem Bi9.1 and $\gamma \in{\{10,30,50\}}$}
\label{Table_exp2_9_1_UB}
\end{table}

\begin{table}[htbp]
\centering
\setlength{\tabcolsep}{2pt}
\begin{tabular}[c]{c|cc|cc|cc}
\multicolumn{ 1}{c}{} & 
\multicolumn{ 6}{c}{$GAP_{P_{sub}}\%$} \\ 
\multicolumn{ 1}{c}{No} & 
\multicolumn{ 2}{c}{$\gamma=10$} & 
\multicolumn{ 2}{c}{$\gamma=30$} & 
\multicolumn{ 2}{c}{$\gamma=50$} \\ 
\hline
1 &  5.58 & 0.74& \textbf{4.25} & \textbf{0.74} & \underline{4.43} & 0.74 \\ 
2 & 3.92 & 1.23 & \textbf{3.48} & \textbf{1.15} & \textbf{3.40} & \textbf{1.13} \\
3 & 1.73 & 2.94 & \textbf{1.59} & \textbf{2.51} & \textbf{1.57} & \textbf{2.42} \\  
4 & 2.11 & 1.51 & \textbf{2.10} & \textbf{1.50} & \textbf{2.01} & \textbf{1.46} \\
5 & 1.78 & 1.58 & \textbf{1.77} & \textbf{1.38} & \textbf{1.68} & \textbf{1.34} \\

\end{tabular}
\caption{Chute2, $GAP_{P_{sub}}\%$ for test problem Bi9.1 and $\gamma \in{\{10,30,50\}}$}
\label{Table_exp2_9_1_GAP}
\end{table}

\begin{table}[htbp]
\centering
\setlength{\tabcolsep}{2pt}
\begin{tabular}[c]{c|c|c|c|c|c|c}
\multicolumn{ 1}{c}{} & 
\multicolumn{ 3}{c}{$|S_U|$} & 
\multicolumn{ 3}{c}{Time$\,S_U (s)$}\\
\multicolumn{ 1}{c}{No} & 
\multicolumn{ 1}{c}{$\gamma=10$} & 
\multicolumn{ 1}{c}{$\gamma=30$} & 
\multicolumn{ 1}{c}{$\gamma=50$} & 
\multicolumn{ 1}{c}{$\gamma=10$} & 
\multicolumn{ 1}{c}{$\gamma=30$} & 
\multicolumn{ 1}{c}{$\gamma=50$}  \\ 
\hline
1 & 11 & 31 & 52 & 439.96 (411.14) & 515.48 (402.94) & 665.55 (404.57)\\  
2 & 10 & 27 & 44 & 482.23 (407.04) & 678.38 (402.56) & 765.19 (409.44)\\   
3 & 8 & 20 & 32 & 443.50 (406.40) & 534.97 (410.38) & 667.23 (404.69)\\
4 & 5 & 15 & 23 & 424.95 (403.06) & 465.13 (400.30) & 484.92 (402.56)\\ 
5 & 5 & 13 & 20 & 432.90 (400.51) & 485.67 (404.07) & 512.03 (403.40)\\ 

\end{tabular}
\caption{Chute2, $|S_U|$, and Time$\,S_U (s)$ for test problem Bi9.1 and $\gamma \in{\{10,30,50\}}$}
\label{Table_exp2_9_1_Su_Time}
\end{table}

\begin{table}[htbp]
\centering
\setlength{\tabcolsep}{2pt}
\begin{tabular}[c]{c|ccc|ccc|ccc}
\multicolumn{ 1}{c}{} & 
\multicolumn{ 9}{c}{$U(S_U,\lambda)$} \\
\multicolumn{ 1}{c}{No} & 
\multicolumn{ 3}{c}{$\gamma=10$} & 
\multicolumn{ 3}{c}{$\gamma=30$} & 
\multicolumn{ 3}{c}{$\gamma=50$} \\ 
\hline
1 & 128872.00 & 131116.00 & 122259.00 & 128872.00 & 131116.00 & 120988.00 & 128872.00 & 131116.00 & 120434.00 \\
2 & 125477.00 & 125196.00 & 121464.00 & 125483.00 & 125196.00 & 120230.00 & 125477.00 & 125196.00 & 119034.00 \\ 
3 & 122652.00 & 131116.00 & 131738.00 & 122094.00 & 131116.00 & 131738.00 & 122334.00 & 131116.00 & 131738.00 \\ 
4 & 122671.00 & 123772.00 & 125154.00 & 122242.00 & 123201.00 & 124765.00 & 122404.00 & 123097.00 & 124677.00 \\ 
5 & 128872.00 & 119155.00 & 131738.00 & 128872.00 & 118229.00 & 131738.00 & 128872.00 & 117837.00 & 131738.00 \\

\end{tabular}
\caption{Chute2, upper bounds for test problem Three6.1 and $\gamma \in{\{10,30,50\}}$}
\label{Table_exp2_3crit_6_1_UB}
\end{table}

\begin{table}[htbp]
\centering
\setlength{\tabcolsep}{2pt}
\begin{tabular}[c]{c|ccc|ccc|ccc}
\multicolumn{ 1}{c}{} & 
\multicolumn{ 9}{c}{$GAP_{P_{sub}}\%$} \\ 
\multicolumn{ 1}{c}{No} & 
\multicolumn{ 3}{c}{$\gamma=10$} & 
\multicolumn{ 3}{c}{$\gamma=30$} & 
\multicolumn{ 3}{c}{$\gamma=50$} \\ 
\hline
1 & 7.95 & 1.91 & 28.07 & 7.95 & 1.91 & \textbf{27.31} & 7.95 & 1.91 & \textbf{26.98} \\
2 & 1.05 & 1.32 & 4.53 & \underline{1.06} & 1.32 & \textbf{3.55} & \textbf{1.05} & 1.32 & \textbf{2.58} \\ 
3 & 26.23 & 2.92 & 7.80 & \textbf{25.89} & 2.92 & 7.80 & \textbf{26.04} & 2.92 & 7.80 \\ 
4 & 1.37 & 1.81 & 1.27 & \textbf{1.03} & \textbf{1.36} & \textbf{0.97} & \underline{1.16} & \textbf{1.27} & \textbf{0.90} \\ 
5 & 10.28 & 9.24 & 1.75 & 10.28 & \textbf{8.53} & 1.75 & 10.28 & \textbf{8.23} & 1.75 \\

\end{tabular}
\caption{Chute2, $GAP_{P_{sub}}\%$ for test problem Three6.1 and $\gamma \in{\{10,30,50\}}$}
\label{Table_exp2_3crit_6_1_GAP}
\end{table}

\begin{table}[htbp]
\centering
\setlength{\tabcolsep}{2pt}
\begin{tabular}[c]{c|c|c|c|c|c|c}
\multicolumn{ 1}{c}{} & 
\multicolumn{ 3}{c}{$|S_U|$} & 
\multicolumn{ 3}{c}{Time$\,S_U (s)$}\\
\multicolumn{ 1}{c}{No} & 
\multicolumn{ 1}{c}{$\gamma=10$} & 
\multicolumn{ 1}{c}{$\gamma=30$} & 
\multicolumn{ 1}{c}{$\gamma=50$} & 
\multicolumn{ 1}{c}{$\gamma=10$} & 
\multicolumn{ 1}{c}{$\gamma=30$} & 
\multicolumn{ 1}{c}{$\gamma=50$}  \\ 
\hline
1 & 8 & 20 & 31 & 106.30 (104.24) & 133.88 (121.32) & 118.75 (102.00)\\ 
2 & 4 & 11 & 17 & 413.24 (408.25) & 436.70 (409.23) & 429.80 (401.94)\\  
3 & 10 & 27 & 44 & 50.87 (47.43) & 76.35 (61.08) & 60.59 (41.18)\\  
4 & 3 & 6 & 10 & 293.45 (289.56) & 380.25 (365.56) & 353.60 (319.88)\\ 
5 & 11 & 30 & 50 & 56.90 (51.22) & 66.93 (51.58) & 79.55 (50.46)\\ 

\end{tabular}
\caption{Chute2, $|S_U|$, and Time$\,S_U (s)$ for test problem Three6.1 and $\gamma \in{\{10,30,50\}}$}
\label{Table_exp2_3crit_6_1_Su_Time}
\end{table}

\begin{table}[htbp]
\centering
\setlength{\tabcolsep}{2pt}
\begin{tabular}[c]{c|ccc|ccc|ccc}
\multicolumn{ 1}{c}{} & 
\multicolumn{ 9}{c}{$U(S_U,\lambda)$} \\
\multicolumn{ 1}{c}{No} & 
\multicolumn{ 3}{c}{$\gamma=10$} & 
\multicolumn{ 3}{c}{$\gamma=30$} & 
\multicolumn{ 3}{c}{$\gamma=50$} \\ 
\hline
1 & 115765.00 & 115804.00 & 113089.00 & 115250.00 & 115546.00 & 112742.00 & 115126.00 & 115335.00 & 112656.00 \\ 
2 & 114092.00 & 113121.00 & 118122.00 & 113842.00 & 112477.00 & 117089.00 & 113861.00 & 112778.00 & 118122.00 \\  
3 & 116240.00 & 117966.00 & 112426.00 & 116229.00 & 117577.00 & 111977.00 & 116160.00 & 117657.00 & 111578.00 \\  
4 & 112893.00 & 116633.00 & 114061.00 & 112291.00 & 116382.00 & 113713.00 & 112177.00 & 116438.00 & 113665.00 \\ 
5 & 119379.88 & 112286.00 & 118122.00 & 119379.88 & 111440.00 & 118122.00 & 119379.88 & 111242.00 & 118122.00 \\ 

\end{tabular}
\caption{Chute2, upper bounds for test problem Three9.1 and $\gamma \in{\{10,30,50\}}$}
\label{Table_exp2_3crit_9_1_UB}
\end{table}

\begin{table}[htbp]
\centering
\setlength{\tabcolsep}{2pt}
\begin{tabular}[c]{c|ccc|ccc|ccc}
\multicolumn{ 1}{c}{} & 
\multicolumn{ 9}{c}{$GAP_{P_{sub}}\%$} \\ 
\multicolumn{ 1}{c}{No} & 
\multicolumn{ 3}{c}{$\gamma=10$} & 
\multicolumn{ 3}{c}{$\gamma=30$} & 
\multicolumn{ 3}{c}{$\gamma=50$} \\ 
\hline
1 & 3.36 & 3.40 & 3.36 & \textbf{2.93} & \textbf{3.19} & \textbf{3.06} & \textbf{2.82} & \textbf{3.01} & \textbf{2.99} \\ 
2 & 3.36 & 8.14 & 3.06 & \textbf{3.14} & \textbf{7.61} & \textbf{2.21} & \underline{3.16} & \underline{7.86} & \underline{3.06} \\  
3 & 4.90 & 3.03 & 4.50 & \textbf{4.89} & \textbf{2.71} & \textbf{4.12} & \textbf{4.83} & \underline{2.78} & \textbf{3.78} \\  
4 & 6.87 & 2.31 & 2.84 & \textbf{6.37} & \textbf{2.10} & \textbf{2.54} & \textbf{6.27} & \underline{2.15} & \textbf{2.50} \\ 
5 & 5.75 & 100.00 & 4.09 & 5.75 & 100.00 & 4.09 & 5.75 & 100.00 & 4.09 \\ 

\end{tabular}
\caption{Chute2, $GAP_{P_{sub}}\%$ for test problem Three9.1 and $\gamma \in{\{10,30,50\}}$}
\label{Table_exp2_3crit_9_1_GAP}
\end{table}

\begin{table}[htbp]
\centering
\setlength{\tabcolsep}{2pt}
\begin{tabular}[c]{c|c|c|c|c|c|c}
\multicolumn{ 1}{c}{} & 
\multicolumn{ 3}{c}{$|S_U|$} & 
\multicolumn{ 3}{c}{Time$\,S_U (s)$}\\
\multicolumn{ 1}{c}{No} & 
\multicolumn{ 1}{c}{$\gamma=10$} & 
\multicolumn{ 1}{c}{$\gamma=30$} & 
\multicolumn{ 1}{c}{$\gamma=50$} & 
\multicolumn{ 1}{c}{$\gamma=10$} & 
\multicolumn{ 1}{c}{$\gamma=30$} & 
\multicolumn{ 1}{c}{$\gamma=50$}  \\ 
\hline
1 & 8 & 20 & 31 & 519.63 (420.77) & 548.70 (400.11) & 654.78 (411.71)\\ 
2 & 12 & 32 & 55 & 412.71 (403.34) & 478.80 (409.01) & 521.99 (401.65)\\ 
3 & 12 & 32 & 52 & 444.12 (400.92) & 527.00 (404.25) & 605.74 (405.75)\\ 
4 & 9 & 22 & 36 & 425.02 (400.11) & 462.81 (405.86) & 515.93 (405.72)\\ 
5 & 5 & 12 & 20 & 325.93 (323.31) & 335.62 (324.23) & 360.17 (330.95)\\ 

\end{tabular}
\caption{Chute2, $|S_U|$, and Time$\,S_U (s)$ for test problem Three9.1 and $\gamma \in{\{10,30,50\}}$}
\label{Table_exp2_3crit_9_1_Su_Time}
\end{table}

\subsection{Comparing Chute2 with Chute1}
\label{comparison}
When comparing Chute2 to Chute1, for all tested instances and all values of the $\gamma$ parameter, we observe no deterioration of any component of $G_{P_{sub}}(R(S_L, S_U,\lambda))$. We observe the following.

For instance Bi6.1, for all values of $\gamma$, all components of $G_{P_{sub}}(R(S_L, S_U,\lambda))$ improve for all vectors $\lambda$.

For instance Bi9.1, for $\gamma=10$, all components of $G_{P_{sub}}(R(S_L, S_U,\lambda))$   improve for four vectors $\lambda$, and one component improves for one vector $\lambda$. The same situation occurs for $\gamma=30$ and $\gamma=50$.

For instance Three6.1, for $\gamma=10$, at least one component of $G_{P_{sub}}(R(S_L, S_U,\lambda))$ improves for all vectors $\lambda$, and for two ones all components improve. The same situation occurs for $\gamma=30$ and $\gamma=50$.

For instance Three9.1, for all $\gamma$, at least one component of $G_{P_{sub}}(R(S_L, S_U,\lambda))$ improves for four vectors $\lambda$. All components of $G_{P_{sub}}(R(S_L, S_U,\lambda))$ improve for $\gamma=10$, $\gamma=30$, and $\gamma=50$, respectively, for  three, four, and three vectors $\lambda$. For all values of the $\gamma$ parameter, for one vector $\lambda$ there is no improvement of any component of $G_{P_{sub}}(R(S_L, S_U,\lambda))$.

Table \ref{Table_AVG_Times} shows times of deriving upper shells averaged over all vectors $\lambda$ for both tested algorithms. We observe a significant increase in these times for  Chute2 compared to Chute1. 
It should be recalled here that Chute2 uses the Suboptimal algorithm, for which the stopping condition depends on the assumed for this algorithm time limit $T^S=400$ seconds.
For Chute2, the average running time of the Suboptimal algorithm is given in parentheses.
It can be seen that in this case a significant fraction of its running time is that of the Subotimal algorithm.
In addition, for all $\gamma$ values, the average running times of Chute2 are larger for Bi9.1 than for Three9.1, which is theoretically a harder problem to solve because it has one more objective function. The implication is that for Bi9.1, for all five lambda vectors, the Subotimal algorithm terminated due to the $T^S$ limit, while for Three9.1 -- for four lambda vectors. This affected the average times. For details, see Tables \ref{Table_exp2_9_1_Su_Time} and \ref{Table_exp2_3crit_9_1_Su_Time}.


\begin{table}[htbp]
\centering
\begin{tabular}[c]{c|c|c}
\multirow{2}{1em}{$\gamma$}  & AVG Time$\,S_U(s)$ & AVG Time$\,S_U(s)$\\
& Chute1 & Chute2\\
\hline
\multicolumn{ 3}{c}{Bi6.1}\\
\hline
10 & 2.09 & 77.47 (76.07)\\
30 & 5.95 & 85.03 (81.60)\\
50 & 10.04 & 83.28 (77.95)\\
\hline
\multicolumn{ 3}{c}{Bi9.1}\\
\hline
10 & 4.02 & 444.71 (405.63)\\
30 & 12.72 & 535.92 (404.05)\\
50 & 21.05 & 618.98 (404.93)\\
\hline
\multicolumn{ 3}{c}{Three6.1}\\
\hline
10 & 4.04 & 184.15 (180.14)\\
30 & 7.98 & 218.82 (201.75)\\
50 & 19.50 & 208.46 (183.09)\\
\hline
\multicolumn{ 3}{c}{Three9.1}\\
\hline
10 & 44.72 & 425.48 (389.69)\\
30 & 172.96 & 470.58 (388.69)\\
50 & 263.05 & 531.72 (391.16)\\
\hline
\end{tabular}
\caption{Average times of deriving upper shells for Chute1 and Chute2.}
\label{Table_AVG_Times}
\end{table}

\subsection{Discussion}
\label{discussion}

For all test instances of the MOMIP problem, with time limits set, algorithm Chute2 determines tighter upper bounds measured with the help of $G_{P_{sub}}(R(S_L, S_U,\lambda))$ than algorithm Chute1 in most cases. Yet, this comes at the expense of a significant increase in the computation time for deriving upper shells. So, we can observe a trade-off between the quality of the interval representation of the implicit Pareto optimal outcome for a given $\lambda$ and computation time.
In both the algorithms, for a given $\lambda$, tightness of upper bounds can be controlled by changing values of parameter $\gamma$. However, changing the $\gamma$ value from lower to higher does not always guarantee an improvement in at least one component of $G_{P_{sub}}(R(S_L, S_U,\lambda))$. It may even happen that some of its components deteriorate. However, in all tested instances, when changing from the lowest to the highest value of parameter $\gamma$, no deterioration of any component of $G_{P_{sub}}(R(S_L, S_U,\lambda))$ has been recorded for all vectors $\lambda$.

The deterioration of some of the components of $G_{P_{sub}}(R(S_L, S_U,\lambda))$ after increasing $\gamma$ may be due to the fact that increasing the value of $\gamma$ does not preserve the elements of the set $S_U$ obtained for smaller $\gamma$, but generates a new, denser set $S_U$, but different in general. These new $S_U$ elements may not be able to generate always better, but in some cases generete even slightly worse vectors of upper bounds than those obtained for smaller $\gamma$.
During decision-making, one can store all the derived upper shells and use their elements in the Chute algorithm as helpers to determine tighter upper bounds for a given $\lambda$ when the DM asks for them.

Parameters affecting the operation of algorithms Chute1 and Chute2 (in particular, time limits for optimization  as well as parameter $\gamma$ ) were arbitrarily set for the numerical experiments conducted on the selected test instances. We can not recommend the adopted parameter values (e.g., $T^L=1200$ seconds, $T^S=400$ seconds) for other instances of the MOMIP problem. The values of these parameters might depend on the problem to be solved, the available computational resources and the conditions of the decision-making process itself.

As Chute1 and Chute2 use a MIP solver as a black box, it is difficult to provide their theoretical performance, especially since they can work with any instance of the MOMIP problem that meets the very generic assumptions made in this work.
During their operation, multiple instances of the single-objective MIP problem are solved, which are parameterized by $\lambda$ in the case of Chute1 and ($\lambda, \mu$) in the case of Chute2.
Moreover, Chute2 uses the Suboptimal algorithm as a black box, and it is difficult to predict which termination condition of Suboptimal will occur as it runs for different instances of the single-objective MIP problem parametrized by $\lambda$.

\begin{figure}
\includegraphics[scale=0.25]{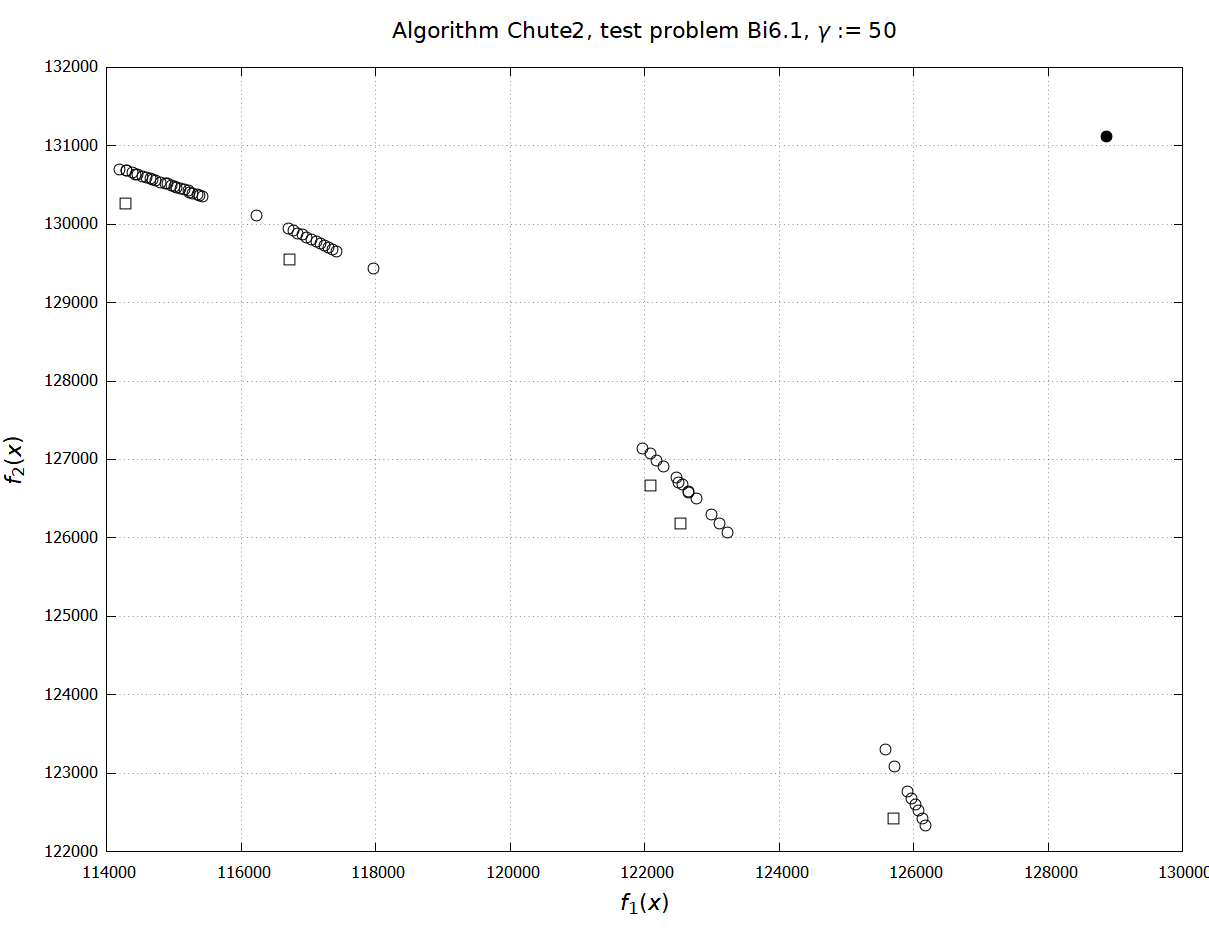} 
\centering
\caption{A finite two-sided approximation of the Pareto front: $\square$, $\circ$ -- images of lower shell $S_L$ and upper shell $S_U$ elements in the objective space, respectively, $\bullet$ -- $y^{*}$.}
\label{Fig_twosided}
\end{figure}

The Chute algorithm returns not only the interval representation but also lower and upper shells. Let us assume that for a given set $\Lambda:=\{\lambda^1,\lambda^2,\lambda^3\}$ the algorithm derives upper shells $S_{U}(\lambda^1)$, $S_{U}(\lambda^2)$, and $S_{U}(\lambda^3)$. 
$S_U:={\widehat{\oplus}_{i=1}^{3}S_{U}(\lambda^i)}$ (where "$\widehat{\oplus}$" is an operator of adding two sets and removing dominating elements) is an upper shell, and $S_L:={\oplus}_{i=1}^{3}\{INC^{\lambda^i}\}$ (where "$\oplus$" is an operator of adding two sets and removing dominated elements) is a lower shell.
One can use $S_L$ and $S_U$ to calculate interval representations of implicit Pareto optimal outcomes designated by $\lambda \notin \Lambda$.
For test problem Bi6.1 and $\gamma := 50$, images of the lower and upper shells obtained this way (for five considered vectors $\lambda$) are shown in Figure \ref{Fig_twosided}. These images form a finite two-sided approximation of the Pareto front.
The approximation does not fully cover the entire Pareto front, as it was derived to determine interval representations of implicit Pareto optimal outcomes designated by just the selected five vectors $\lambda$.

We can say that we obtained the two-sided approximation of the Pareto front, shaped by the DM's preferences expressed with the help of vectors $\lambda$.

Although we aim not to derive
approximations of the entire Pareto front 
(as in multi-objective branch and bound, see, e.g., \cite{Przybylski_Gandibleux_2017}, \cite{Forget_et_all_2022}), 
with a fairly large set of evenly distributed  vectors $\lambda$, one
would imagine (for $k=2,3$) the corridor in which the Pareto front is located.
\section{Limitations of the Chute algorithm and its possible enhancements}
\label{limitations}
In the Chute algorithm, we have assumed that for all probing vectors $\lambda^{'}$ in the FindUpperShell algorithm, the same vector of multipliers $\mu$ is used. The Chute1 version inherently uses a single vector $\mu$. Yet, for the Chute2 version, it is just a heuristic assumption that vector $\mu$, set with the help of the Suboptimal algorithm for a given vector $\lambda$ in line 4 of the Chute algorithm, provides a tight lower bound on values of the objective function of problem (\ref{eq2.2}) for probing vectors $\lambda^{'}$ close to $\lambda$. However, this need not be the case, especially for vectors $\lambda^{'}$, which are significantly different from $\lambda$ (i.e., when they indicate a significantly different search direction in the objective space).

However, one can imagine version Chute3 of the Chute algorithm in which the determination of vector $\mu$ takes place in the FindUpperShell algorithm for each probing $\lambda^{'}$ considered in it  (or, e.g., for $\lambda^{'}$ not close, in a sense, to a given $\lambda$).  This, at the same time, would require adopting a reasonable time limit on optimization in the Suboptimal algorithm, as we expect many probing vectors $\lambda^{'}$ in the FindUpperShell algorithm. This time limit could be, e.g., a fraction of time $T^S$ adopted in Chute2.
Since the number of probing vectors $\lambda^{'}$ is not a priori known, this time limit would have to be determined by some heuristic rule. It is not desirable that excessive time to determine all vectors $\mu$ is a barrier to the applicability of the proposed method.

In a real decision-making process using the Chute algorithm, it is possible to calculate a more adjusted value of parameter $\gamma$ for a new vector $\lambda$ based on the properties of the lower and upper shells obtained for previous vectors $\lambda$, and with which this algorithm was called. Let us look, for example, at Table \ref{Table_exp2_6_1_Su_Time}. For $\lambda^{1} = (0.055,\,0.945)$, $|S_U|=26$, Time$\,S_U = 91.44$ seconds, but for $\lambda^{5} = (0.439,\,0.561)$, $|S_U|=7$, Time$\,S_U = 39.82$ seconds. To have the time of deriving an upper shell for some  $\lambda^{'}$ close to $\lambda^{1}$ comparable to the time for $\lambda^{5}$, it could be possible to lower the value of parameter $\gamma$ from 50 to, e.g., 20. Such an overarching mechanism (with a set of rules based on statistics collected during the decision-making process) for controlling the behavior of the Chute algorithm could be useful when computation time is an important factor.

In the proposed method of deriving upper shells (algorithm FindUpperShell), there is no single parameter to limit optimization time for getting theirs. Yet, such a time limit can be incorporated relatively easily as an additional stop condition in FindUpperShell.
More generally, it could be even desirable to introduce in the Chute algorithm a time limit for determining the interval representation of the implicit Pareto optimal outcome for a given vector $\lambda$.
The DM would give, for example, in addition to $\lambda$ and $T^L$, time limit $T^I$ (e.g., $T^L=1200$ seconds, $T^I=500$ seconds). Then the Chute algorithm would have time limit $\frac{T^I}{k}$ to derive an upper shell for calculating a single component of the interval representation of implicit Pareto optimal outcome designated by $\lambda$.  Determination of suboptimal vectors $\mu$ in Chute2 and Chute3 would, of course, have to be within some fraction of $T^I$.

Based on the above alone, one can imagine many schemes for budgeting calculations, leading to providing interval representations in a decision-making system based on the Chute algorithm.

In our approach, to find elements of the upper shell we solve to optimality the Chebyshev scalarization of the surrogate relaxation of the MOMIP problem. For instances of the MOMIP problem with a large number of constraints (e.g., 1000), even with a suboptimal vector of multipliers $\mu$ provided by the Suboptimal algorithm (that is, with a single constraint that mimics the original set of constraints of the MOMIP problem), the FindUpperShell procedure may not derive elements $x$ of the upper shell that $f_l(x) < y^{*}_l$, $l=1,\ldots,k$.  In this case, the upper bounds on components of Pareto optimal outcome designated by $\lambda$ are not better than components of $y^{*}_l$. That is, images of elements of the upper shell in the objective space are very far from the Pareto front of the MOMIP problem, and do not provide better upper bounds than the components of $y^{*}$.

To find (sub)optimal values of multipliers $\mu_p$, other algorithms can be used (see, e.g., \cite{Sikorski_1986}). To find elements of upper shells, sophisticated combined relaxation techniques for MIP problems, e.g., Lagrangean/surrogate heuristics (see \cite{Narciso_Lorena_1999}) can also be applied. In the current work, we consider the most general formulation of the MOMIP problem, but to find those elements, problem-specific techniques may help. 
The disadvantage of the proposed generic scheme Chute is that it does not take into account the specifics of a given instance of the MOMIP problem. However, by showing its Chute2 modification and pointing to the Chute3 option, it has been shown how this scheme can be modified.

Within the generic framework presented, other methods of deriving upper shells in the FindUpperShell procedure can also be applied, e.g., a method shown in \cite{Miroforidis_2021_Markowitz}.

\section{Final remarks}
\label{final_remarks}
It has been shown how to algorithmically derive lower and upper shells to the MOMIP problem (for any $k>1$) to get the interval representation of the Pareto optimal outcome designated by vector $\lambda$. 
On selected examples, it has been shown that with the help of the proposed method, one can find such interval representations for randomly selected vectors $\lambda$ where there is a time limit for a MIP solver on deriving a single Pareto optimal solution.

We conducted some preliminary experiments with the Chute algorithm on instances of the MOMKP with four objective functions. 
However, due to the mechanism adopted in the FindUpperShell algorithm for changing the probing $\lambda$ vectors, the results achieved were not satisfactory.



In our future work, we want to improve the method of populating an upper shell (in quest of finding its elements that can provides upper bounds) by changing the scheme of probing the objective space.
We want it to determine upper shells with the desired properties for four and more objective functions.
We also want to apply the presented generic approach to other instances of the MOMIP problem, especially ones connected to real-life problems. This would help verify the practicality of the proposed general method and identify those elements that could be tailored for specific instances of this problem.
Possible modifications to the proposed method are indicated in Section \ref{limitations}. These are also worth considering in further work.



\bibliographystyle{infor}
\bibliography{biblio}
\vfill


\begin{biography}\label{bio1}
\author{G. Filcek}  is an assistant professor at the Wroc\l{}aw University of Science and Technology, Poland. He obtained his Ph.D. in Computer Science from that institution in 2011. His scientific interests include multi-objective optimization, multiple criteria decision-making, intelligent decision support systems, blockchain technologies, and systems engineering.
\end{biography}

\begin{biography}\label{bio2}
\author{J. Miroforidis} is an assistant professor at the Systems Research Institute, Polish Academy of Sciences, Poland.
He obtained his Ph.D. in Computer Science from that institution in 2010.
His scientific interests are in multi-objective optimization, multiple criteria decision making, and computer-aided
decision making.
\end{biography}

\end{document}